\documentclass{amsart}
\usepackage{graphics,verbatim,color,amsthm}

\usepackage[all]{xy} 

\usepackage{hyperref}

\theoremstyle{plain}
\newtheorem{theorem}{Theorem}
\newtheorem{proposition}{Proposition}

\newtheorem{lemma}{Lemma}  
\newtheorem{cor}{Corollary} 
\newtheorem{definition}{Definition}

\theoremstyle{definition}
\newtheorem{example}{Example}

\newtheorem{exercise}{Exercise}

\usepackage{enumitem}

\newlist{noitemize}{itemize}{1}
\setlist[noitemize]{label={}, labelsep=0pt, leftmargin=0pt}

\def\R{\mathbb{R}}
\def\Z{\mathbb{Z}}
\def\C{\mathbb{C}}

\def\E{\mathbb{E}}
\def\V{\mathbb{V}}
\def\A{\mathbb{A}}

\def\eps{\epsilon}

\def\bS{\mathbb{S}}

\def\Ri{ Riemannian }

 \newcommand{\dd}[2]
{
{{\partial #1}   \over {\partial #2}}
}

\newcommand{\beq}{\begin{equation}}
\newcommand{\eeq}{\end{equation}}

\newcommand{\Leq}[1]{\label{#1}\end{equation}}

\begin{document}

\author{Richard Montgomery}
\address{Mathematics Department\\ University of California, Santa Cruz\\
Santa Cruz CA 95064}
\email{rmont@ucsc.edu}

\date{September 15, 2016}

\title{The Kepler Cone, Maclaurin Duality  and Jacobi-Maupertuis metrics. }

\begin{abstract}    The Kepler problem is
the special case $\alpha = 1$ of the  power law problem:
to solve   Newton's equations for a   central force    whose 
potential is of the form $-\mu/r^{\alpha}$ where $\mu$ is a coupling constant.  Associated to such a problem
is a two-dimensional cone with cone angle  $2 \pi c$ with  $c = 1 - \frac{\alpha}{2}$. 
We construct a transformation taking  the  geodesics
of this cone   to the  zero energy  solutions  of the $\alpha$-power law problem. 
 The `Kepler Cone'  is the cone  associated to the Kepler problem.  
 This zero-energy cone transformation is a special case of a  
 transformation  discovered by Maclaurin in the 1740s
 transforming  the $\alpha$- power law problem for any energies to a `Maclaurin dual'
 $\gamma$-power law problem where  $\gamma = \frac{2 \alpha}{2-\alpha}$
 and which, in the process,  mixes up the energy of one problem with the coupling constant of the other.
 We derive Maclaurin duality using
 the  Jacobi-Maupertuis metric reformulation of mechanics. 
We then use the conical metric   to explain properties of  Rutherford-type scattering off power law
potentials at positive energies.  The one   possibly  new result in the paper
concerns  ``star-burst curves''  which arise as   limits of families 
negative energy solutions as their angular momentum tends to zero.
We also describe some  history around Maclaurin duality and
 give two derivations of the   Jacobi-Maupertuis
metric reformulation of classical mechanics.   The piece is expository,
aimed at an upper-division undergraduate. Think  American Math. Monthly. 
\end{abstract} 

\maketitle

\section{The Kepler Cone}
\begin{figure}[h]
\scalebox{0.4}{\includegraphics{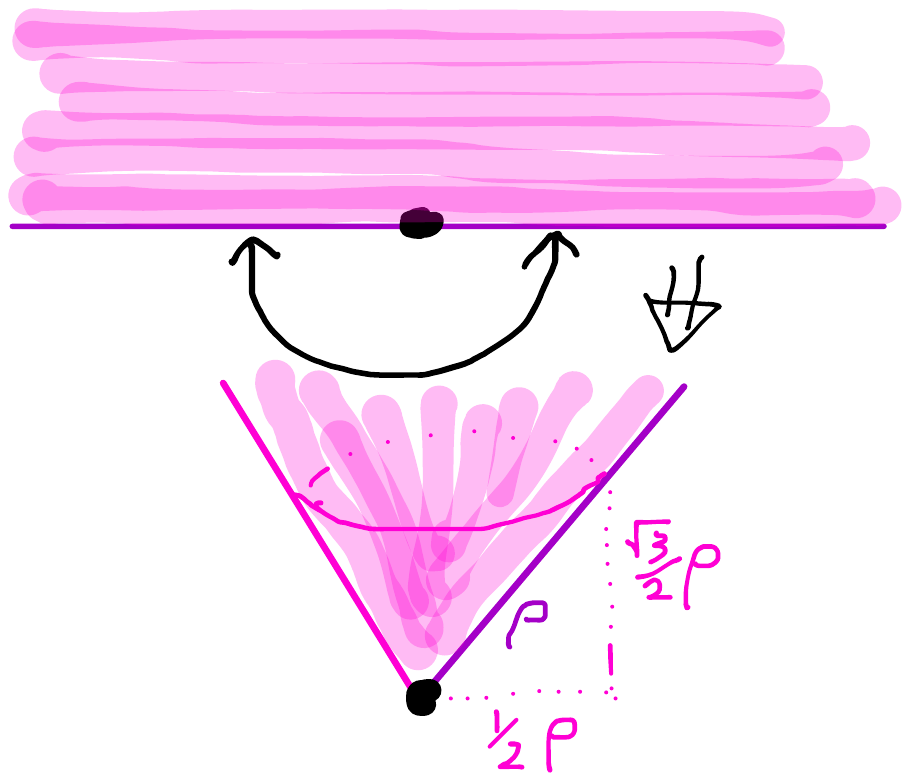}}
\caption{Folding a sheet of paper to make a cone.}  \label{cone1}
\end{figure}

Take a   sheet of paper.  Join   adjacent corners, creasing and folding so as to  bisect  their common edge. 
Join the resulting two half-edges together and tape   together.    You have
formed part of the Kepler cone with the old edge forming one generator.  Lines drawn on the
paper    become geodesics on the cone.  See figure \ref{cone1}. 
 The cone  embeds  isometrically  in 3-space as  the locus 
$$x^2 + y^2 = 3 z^2,  z \ge 0. $$

Place the Kepler cone  off  to the side.  Take out a new sheet of paper
representing a Euclidean plane which  we   call the Newtonian plane.
Mark an origin on this plane.    Draw some parabolas  on the plane with   origin
as   focus.     When parameterized according 
to Kepler's second law  (``equal times in equal areas'') these  parabolas are  
the zero-energy solutions to Kepler's problem
\beq
\ddot q = - \mu \frac{q}{r^3},  r = \| q \|: = \sqrt{ x^2 + y^2}. 
\Leq{Kepler} 
Here  $\mu$ is a positive  constant  called the ``coupling constant''. 
The vector $q = (x,y) \in \R^2$   coordinatizes the Newtonian plane. 
A solution to Kepler's problem is a curve $t \mapsto q(t)$ parameterized by Newtonian time $t$.
 Dots over $q$
denote time derivatives   so that   in the equation  $\ddot q$  denotes the  acceleration
or  second derivative. And $\dot q$ denotes the  velocity or first derivative along  the curve.   The  energy 
associated to Kepler's problem is 
 \beq
 E = \frac{1}{2} | \dot q|^2 - \frac{\mu}{r}.
 \Leq{Kepler E}
and   is constant along solutions to equation (\ref{Kepler}).   
The Kepler parabolas, and their degenerations, the rays through the origin,
exhaust  the supply of zero energy solutions to Kepler's problem.  

Retrieve the Kepler cone.  Consider
the plane  containing the sheet of paper  used to make
the cone.  Call this   the   folding plane.  
Put Cartesian coordinates $Q= (u,v)$ on the folding plane
so that the origin $O$  is the midpoint of the edge of the sheet paper and the sheet lies in the
half-plane $v \ge 0$.     When we fold and glue we identify
  points along the bounding line $v=0$ by the    isometry 
$Q \mapsto -Q$.  The   two rays    merge to form  a single  generator of the Kepler cone  which we will call  {\it the seam}. Rotate the half-plane and we  get another representation of the Kepler cone,
one with a different seam.  When we use this new half-plane   the gluing map still glues  
$Q$ to $-Q$ on its bounding line.  If $Q$ is in the interior of the half-plane  then $-Q$ is in
the complement of the closed half-plane.  We have derived   the following   algebraic construction of the Kepler cone.
\begin{lemma}
The Kepler cone is the metric quotient   $\R^2 / (\pm 1)$ of $\R^2$
by the action of the two element group $\pm 1$ acting by $Q \mapsto \pm Q$.
\end{lemma}
\noindent We describe in subsection \ref{subsec: metric quotient} what we mean by  ``metric quotients'', that is to say,  how the Euclidean
 metric on $\R^2$ induces a metric on the quotient.

Identify   the folding plane and the Newtonian plane with a copy  of 
$\C$ by writing $q= (x,y) = x+iy$ and $Q = (u,v) = u+ iv$.  Then 
 squaring 
 \beq
Q \mapsto q = Q^2,
\Leq{folding map}
 defines a map from the folding plane to the Newtonian plane.   
  Because    $(-Q)^2 = Q^2$, the  squaring  map
    induces a map from  the  Kepler cone  to the Kepler plane, namely the map   $\{Q, -Q \} \to Q^2$.  
 
  \begin{exercise}  Show that squaring (\ref{folding map}) takes 
straight lines   to   Kepler parabolas.
\end{exercise}

\noindent Did foisting this exercise on you   feel  like pulling  a rabbit out of a hat?
Sure, you can do the algebra, but why should it be true?

 Introduce the ``Jacobi-Maupertuis''  metric 
$$ds^2_{1}  = \frac{4 |dq|^2}{|q|} \text{ where } |dq|^2 = dx^2 + dy^2 $$
on the Kepler plane. 
This metric is Riemannian away from $q = 0$ and
the distance function which it induces on the
Newtonian plane is that of a complete metric space.  
\begin{proposition}

 A.  The Kepler parabolas   are   geodesics for the  metric $d^2 s_1$.
 
B. Squaring  (map (\ref{folding map})) induces an isometry between  the Kepler cone and  the Kepler plane endowed with the
Jacobi-Maupertuis metric  $ds^2_{1}$.
\end{proposition}

\noindent Recall that the straight lines in the folding plane became  geodesics on the cone.
Since  isometries  take geodesics to geodesics,    the proposition supplies
a metric  explanation of  our rabbit  out of the hat trick.

Proposition  1 generalizes  to other cones and other force laws.  
A choice of smooth function $V: \R^2 \to \R$  defines a Newton's equations:  
\beq
\ddot q = - \nabla V (q).   
\Leq{Newtons}
These  equations admit a conserved energy  
 $$E(q, v)  = K (v) + V(q), \qquad K(v) = \frac{1}{2} |v|^2,  v = \dot q $$
 as did Kepler's problem whose potential is $V = - \mu/ r$.   

By a central force problem we mean  Newton's equations
 (\ref{Newtons}) for  a rotationally symmetric potential:  
 $V(q) =  f(r)$ where $r = \|q\|$ for some  $f: \R_+ \to \R$.  We are 
 concerned here  with  central force problems of the
 form $V(q) =  -\mu/ r^{\alpha}$.  We call these
 {\it   power law potentials} and the corresponding Newton's equations
 (equation (\ref{eq: C force 1}) below) the {\it power law problem} or 
 {\it power law dynamics}.

\begin{figure}[h]
\scalebox{0.25}{\includegraphics{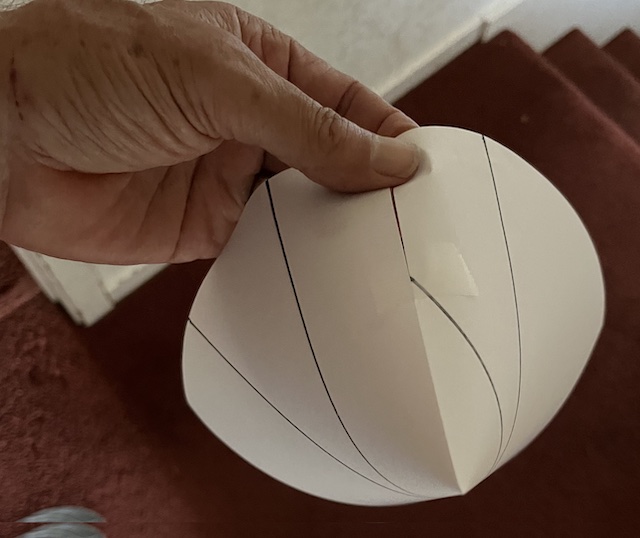}}
\scalebox{0.06}{\includegraphics{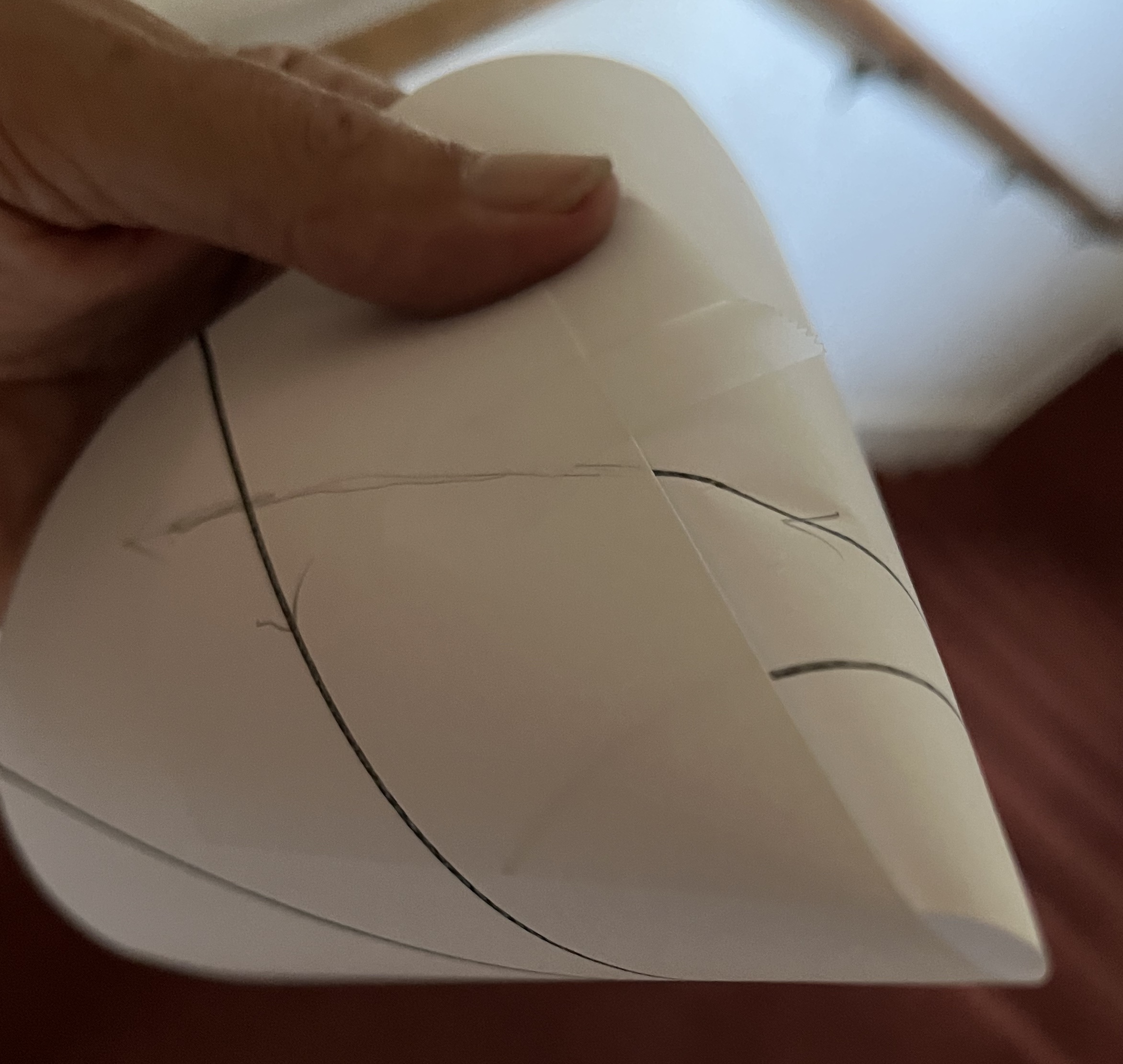}}
\caption{Sharpening the cone angle.}  \label{cones}
\end{figure}

\begin{theorem} \label{thm 1}

 A.  The   geodesics for the JM [Jacobi-Maupertuis]   metric
 \beq
ds^2 _{\alpha}   = \frac{ |dq|^2}{|q|^{\alpha}} 
\Leq{metricalpha}
 are the zero energy solutions to  the power law problem (equation (\ref{eq: C force 1})).
 
 B.     The JM metric   (\ref{metricalpha})  is isometric to the conical metric
 \beq ds^2 _{cone}  = d \rho^2 + c^2 \rho ^2 d \theta ^2,  \text{ with }   c = 1 - \frac{\alpha}{2}
 \Leq{cone metric}    provided   $\alpha \ne 2$.  
 This is the metric for a   two-dimensional cone of angle $|c|$. Here $(\rho, \theta)$ are polar coordinates
 on the plane and are  related to the   polar coordinates $(r, \theta)$ of $q = r e^{i \theta}$
of  the Newtonian plane by $\rho = r^c$.
 
C. Again with  $\alpha \ne 2$,   the  map $Q \mapsto q = Q^{1/c}$ maps straight lines on the folding plane and hence 
geodesics on the  cone of part B,  to the 
zero energy  solutions to Newton's equations described in  part A.
 
\end{theorem}

{\sc Remark on Part B}  The coordinates $(\rho, \theta)$ of part B are
related to  the folding coordinate  $Q \in \C$ of part C 
as follows.  Write $Q  = \rho e^{i \psi}$.  Then  $\theta = \frac{1}{c} \psi$.
For example, for the Kepler cone   $\theta = 2 \psi$ 
and the metric of B is $d \rho^2 + \frac{1}{4} \rho^2 d \theta ^2$.  
See section \ref{sec: cones} for more on cones and for a proof of part B. 

{\sc Remark on $\alpha =2$.}  When $\alpha =2$  the JM metric
is isometric to  that of a cylinder, not a cone.  See subsection \ref{subsec: cylinder} below.

{\sc On folding these cones} We can make all the cones of the
theorem having   $1 \le \alpha < 2$, so $0 < c \le 1/2$, from the same  sheet of paper out
of which we made the Kepler cone. 
Take the two half-edges that we had earlier  glued as soon as  they  had touched
and   continue to    wrap them around tighter so as to 
make   the cone tighter and   sharper   
before taping the half-edges down.  Lines drawn on the paper still make geodesics
on the cone.   Sharper cones means larger coefficients $A^2$ in
the equation of the cone.  See figure \ref{cones}.   

Alternatively,  here is a paper folding construction that
works for all $0 < \alpha < 2$. Cut a sector out of the folding plane,
whose opening angle is $2 \pi (1-c)$ where $c = 1/\beta$,
leaving a sector whose opening angle is $2 \pi c$.  Glue the two
rays bounding this remaining sector together to make any of the cones of
part B of the proposition. (We discuss cones in detail further on.)  

 \begin{figure}
\scalebox{0.4}{\includegraphics{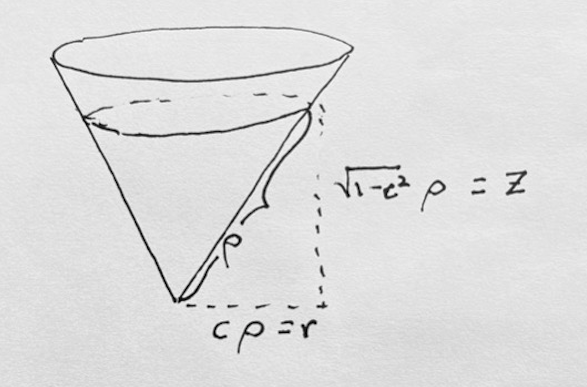}}
\caption{The metric on the cone of part B can be emdedded
as a cone of revolution in three-space provided $0 < c <1$.}  \label{conetilingC}
\label{fig: cone embedded}
\end{figure}

These cones embed as cones of revolution in $\R^3$ provided   $0 < c < 1$.   
See figure \ref{fig: cone embedded}.  Some algebra shows that
the cone is described as
$$x ^2 + y^2 =  A z^2,  z \ge 0 \text{ with }  A = \frac{ 1- c^2}{c^2} .$$
Here $\rho^2 = x^2 + y^2 + z^2$ and $r^2 = x^2 + y^2$.

\vskip .4cm

The Jacobi-Maupertuis principle, presented in the next section,
provides the  engine behind Theorem 1 and its corollary,
Proposition 1.    

\section{Jacobi-Maupertuis}  
\label{sec: JM}

Return to  Newton's equations (\ref{Newtons}).
Allow the potential   $V$ to have poles, i.e. to  take on the values $\pm \infty$.
We will also call  the poles ``collision points'',
in honor of the Kepler problem where $V(q)  = -\infty \iff q = 0$.   In the Kepler problem $q$
records  the planet's location and the origin represents the sun  location so that $q=0$ 
means that   the planet has collided with the sun.  
We assume   $V$ to be continuous everywhere  and    smooth away from collision points.

Recall that the   total  energy 
 $E(q, v)  = K (v) + V(q)$
 is conserved. 
From $K(v)  \ge 0$ we see that   if a   solution $q(t)$ has energy $E$  then along the solution 
 $V(q(t)) \le E$.
 In other words,  solutions having energy $E$ are   constrained to lie
 in the   region $\{q \in \E:  E \ge  V(q)  \}$ within the Newtonian plane.
 We call this region  the ``Hill region''
 for the choice of energy $E$ and its boundary
 $\{q \in \R^2:  E = V(q) \}$ we call  the Hill boundary.

 \begin{definition} The Jacobi-Maupertuis [JM for short]  metric at energy $E$
 for the Newton's equation (\ref{Newtons}) is the  Riemannian metric 
 \beq
 ds^2 _{JM} = 2(E-V(q)) |dq|^2
 \Leq{JM metric} defined on the Hill region $\{ V \le E \}$ and degenerating
 at its boundary $\{ V = E \}$.    
 \end{definition} 
 
  We say that $q_0 = q(t_0)$ is a brake point of a solution $q(t)$  if $\dot q (t_0) = 0$. 
  The solution has ``braked'' to a stop. 
If the energy of this  solution $E$ then being a brake point
is equivalent to hitting the Hill boundary since then  $V(q) = E$   if and only if $K(\dot q (t)) =0$. 
 \begin{theorem} [Jacobi-Maupertuis principle]  Away from   collisions  and brake  points, every energy $E$  solution to Newton's equations
  is the reparameterization of a  geodesic for the energy $E$ JM metric (\ref{JM metric}) on the   Hill region.     Conversely, away from collisions  and   the Hill boundary, every    JM    geodesic is a   reparameterization of 
  an energy $E$ solution of Newton's equations.
 \end{theorem} 
 
 See section \ref{sec: JM proofs}  for   proofs of the theorem.

 {\sc Proof of part A of proposition 1. } Newton  showed
  that the non-collision solutions to Kepler's problem
 satisfy Kepler's 1st law: they are conic sections with one focus at the origin.
 He also showed that the zero energy solutions are the parabolas.
Part A of the proposition now follows immediately from the theorem. 
 
\section{Squaring and Levi-Civita} 

Take the JM metric for Kepler's problem at any energy:
$$ds^2 = 2 (E + \frac{\mu}{r} ) |dq|^2$$
The squaring substitution $q = Q^2$ yields $r = |q| = |Q|^2$,
$dq = 2 Q dQ$ and $|dq|^2 = 4 |Q|^2 |dQ|^2$.
Then $\frac{1}{r} |dq|^2 = 4 |dQ|^2$.  The singularity of the metric at the origin   $r  = 0$ has cancelled out!
Rewritten in terms of $Q$ we find
$$ds^2 = 2( 4 E |Q|^2 + 4 \mu) |dQ|^2$$
We recognize this to be the  JM metric at energy $4\mu$ associated with the potential energy
$W(Q) = -4 E |Q|^2$ on the folding plane.   In going from $V = - \mu/r$ to $W$
the roles of the energy  $E$ and coupling constant $\mu$ have switched!  
 The corresponding
Newton's equations on the folding plane are the linear equations
\beq
Q'' =  8 E Q  
\Leq{Osc}
where the  second derivative $Q''$ is with respect to a `Newtonian time'' $\tau$ on the folding plane.
It follows that the squaring transformation maps solutions of  the linear equation
(\ref{Osc})  to solutions to Kepler's equations having energy $E$. 

{\sc Proof of part B of   proposition 1.}  When  $E = 0$ we have that $ds^2 = 4 \mu |dQ|^2$, which is the flat metric on the
folding plane. Upon folding this metric in turn induces the   metric on the Kepler cone upon folding. The linear ODE
(\ref{Osc})  is that of a straight line
$Q'' = 0$ which fold to geodesics on the cone.    
\vskip .3cm
{\sc Remarks on scaling.}   Multiplying a metric by a positive constant does not change its  
geodesics, so we could just as well use $|dQ|^2$ instead of $4 \mu |dQ|^2$
and correspondingly dilated the Kepler cone metric.  
  Cones, like the Kepler cone,
admits dilations.    However
we scale the metric on the Kepler cone , the circle of radius $\rho$ about the cone point
will have circumference $\pi \rho$, not $2 \pi \rho$.  Scaling the folding plane
by $\lambda$ effects scaling its quotient,  the Kepler cone,  by $\lambda$
\vskip .3cm

We return to the question of the time  $\tau$
used in taking the second derivative defining the ODE   (\ref{Osc})  in the folding plane,
deduced  now  for general energies $E$ in the Newtonian plane.
What is the relation between Newtonian time $t$ and this other time $\tau$? 
If $s$ is JM arclength then we have 
$$\frac{d \tau} {dt}  = \frac{ d \tau} {ds}  \frac {ds } {d t}.$$
Recall   the relation $ds = \lambda dt$ between JM arc-length
$s$ and Newtonian time $t$ when we use the JM principle and where $\lambda = 2(E + V(q))$
is the conformal factor relating the Euclidean metric and the JM metric.    On the $q$ plane with its time $t$
we have $\lambda = 2(E + \mu/ r)$ so  $\frac {ds } {d t} = 2 (E + \mu/r)$.  
 But on the $Q$ plane with its time $\tau$ and   potential  
$W(Q) = -4E|Q|^2 = - 4 E r$ we have   $ds = \tilde \lambda d \tau$
with $\tilde \lambda = 2(4E r + 4\mu) = 4r (2E + \mu/r)$ and so $\frac{d \tau}{ds} = \frac{1}{4r(2 E + \mu/r)}$. 
It follows that $d \tau/ dt = 1/4r$ or  
\beq
d \tau = \frac{dt}{4|q|}
\Leq{LeviCivita2}

{\it  The Levi-Civita transformation} is  the  squaring map (\ref{folding map}) 
together with the time reparameterization (\ref{LeviCivita2}).
See section \ref{sec: history} regarding the long history of an old generalization of this transformation and its name.
  We have shown that
this transformation has the remarkable property that it takes
solutions to Kepler's problem to  solution to the linear ``oscillator'' equation (\ref{Osc}).

 \subsection{Maclaurin  Duality for power law problems.}   
  The Levi-Civita transformation exhibits a duality
 between  the Kepler force law generated by the   $1/r$ gravitational  potential and the   linear force law of Hooke   generated by  an $r^2$ potential.   Maclaurin   knew of this duality centuries before Levi-Civita,
 having had  discovered a    generalization  which applies   to any power law problem.
 The $\alpha$ power law problem, namely  Newton's equations  for the 
 power law potential $V = -\mu r^{-\alpha}$,   is
 \beq 
\ddot q = - \mu \alpha \frac{q}{|q|^{\alpha +2}}.
 \Leq{eq: C force 1} 
 The associated  JM metric at energy  $E$ for this  $V$ is  
$$ds^2 = 2(E + \mu /r^{\alpha})|dq|^2.$$    

The substitution  $q = Q^{\beta}$ yields 
$dq = \beta Q^{\beta-1} dQ$ and $r= |Q|^{\beta}$. The trick is to  tune
the exponent  $\beta$ so  that   the   factor 
$|Q|^{2(\beta -1)}$ occuring in the identity  $|dq|^2 = \beta^2 |Q|^{2(\beta-1)} |dQ|^2$ 
cancels   with  the factor  $\frac{1}{r^{\alpha}}$ arising from the singularity of the potential.   
Thus we require $2(\beta -1) = \beta \alpha$ or  
\beq
 2 - \frac{2}{\beta} = \alpha \text{ or } \beta = \frac{2}{2 - \alpha}
\label{folding exponent}
\eeq 
With  this choice of $\beta$ we find that,   written out in the $Q$-variable,  
\beq
ds^2 = 2(\beta^2 E  |Q|^{\gamma} + \mu \beta^2) |dQ|^2
\label{Q JM}
\eeq
 where 
 \beq
 \gamma = 2 \beta - 2
 = \frac{2 \alpha}{2-\alpha}
 \Leq{dual exponent}
We   recognize expression  (\ref{Q JM})   as    the  JM metric {\it in the folding plane}
for    energy $\mu \beta^2$
associated to the power law potential $W(Q) = - \beta^2 E  |Q|^{\gamma}$.  
Its corresponding Newton's equation is 
  \beq
Q'' = \frac{E \gamma }{\mu \alpha} |Q|^{\gamma -2} Q.
\label{eq: C force 2}
\eeq  
 This  transformed ODE  is  the ``Maclaurin dual'' to our original
 power law problem, equation (\ref{eq: C force 1}).  
 
The transformed   ODE (\ref{eq: C force 2}) is a differential equation  for   curves $Q = Q(\tau)$ in
the folding plane. The new time $\tau$ of $Q'' = \frac{d^2 Q}{d \tau^2}$ 
is not the old Newtonian time $t$ used for the $q$-equation.  The same trick we used when we derived the
reparameterization part of the   Levi-Civita transformation    works here   to find the relation
between  the  $t$-time of the Newtonian plane  and the $\tau$-time of the folding plane.  
This trick used the JM arclength $s$ as an intermediary between the two times.
We have
\begin{eqnarray}
 \frac{d \tau} {dt}  & = & \frac{ d \tau} {ds}  \frac {ds } {d t} \\
 & = &  (2(E \beta^2 |Q|^{\gamma} + \mu \beta^2))^{-1} 2 (E + \mu/r^{\alpha})
 \end{eqnarray}
 Now factor out   $\beta^2 |Q|^{\gamma}$ from
 the first denominator  $2(E \beta^2 |Q|^{\gamma} + \mu \beta^2)$,
   recalling that we tuned $\beta$ so that  $|Q|^{-\gamma} = r^{-\alpha}$.
   We get that   $2(E \beta^2 |Q|^{\gamma} + \mu \beta^2) = \beta^2 r^{\alpha} (2(E + \mu/r^{\alpha})$.
   It follows that   $$\frac{d \tau} {dt} = \frac{1}{ \beta^2 r^{\alpha}}$$
   
 We summarize.
Let the exponents $\alpha, \beta, \gamma$ be related by the relations
(\ref{folding exponent}, \ref{dual exponent}).  
We could equivalently summarize the relations by
 \beq
(1- \frac{\alpha}{2})(1 + \frac{\gamma}{2}) =1, \qquad \frac{\alpha}{2} + \frac{1}{\beta} = 1
\Leq{reln between power laws} 
 Then the    transformation  
  \beq q = Q^{\beta} ,   
   dt = \beta ^2 |q|^{\alpha}  d \tau
  \Leq{Maclaurin   transformation}  
  takes solutions  $q(t)$ to (\ref{eq: C force 2} ) having energy $E$
  to  solutions $Q(\tau)$ to  the   ODE (\ref{eq: C force 2} ).     
   We will call the transformation
  ( \ref{Maclaurin   transformation}) the   ``Maclaurin   transformation'' in honor of Maclaurin  
  who discovered this duality between central force laws whose
  exponents $-\alpha, + \gamma$ are related by 
  (\ref{reln between power laws}).   
  See \cite{60}  section 451 for the statement and   section 875 for the proof.

{\sc Multivalued maps.}  
 
$Q \mapsto Q^{\beta}$ is 
  a  multivalued map  when  $\beta$ is not an integer,  so some care  must be taken in applying the transformation.     We 
  can proceed as follows.  Let  $Q = \rho e^{i \psi}$. 
  Then the  values taken by  $Q^{\beta}$ are any one  of  the possibly countably
  many values  $q = \rho^{\beta} e ^{ i  (\beta \psi + 2 \pi \beta k)}$,
  $k \in \Z$.    Let  $Q(\tau), \tau \in \R$ be an analytic
  curve in the plane which misses the origin.   Choose a point, say $Q(0)$
  along this  curve and one of these  values,   $q_0$ for $Q(0)^{\beta}$.   Analytically continue $Q^{\beta}$ along 
  the entire curve.    In this way we construct an image  curve $q(\tau) = Q^{\beta} (\tau)$ passing through $q_0$.
 Reparameterize $q(\tau)$  by $t$ by doing the integral defining the relation between $t$ and $\tau$
  so as to get the Maclaurin transformed curve $q(t)$.   Different choices of initial value for 
  $q_0$ will be related by rotation by some angle $2 \pi \beta k$.  The corresponding   
  transformed curve through this new $q_0$    will   be related to the original curve  by this same rotation.
  It is still a solution since rotations act as symmetries of central force problems. 
    
The inverse process is identical.  Start with a  curve  $q(t)$ and choose a branch of $q \mapsto q^{1/\beta} =Q$
  near $q(0)$ 
  and   analytically continue $q^{1/\beta}$  along $q(t)$.  We have a unique analytic continuation  as long as the curve
  $q(t)$  misses the origin.  Then we can reparameterize it to get a curve  $Q(\tau)$.      In this way we have a
  correspondence between the non-collision solutions to the two ODEs.
  
  It is worth remarking that the collision solutions, being rays through the origin  also  correspond 
  to each other under the Maclaurin   transformation.

\section{An Abbreviated  History of Maclaurin's  duality. }
\label{sec: history}

In an early  draft of this article I used   ``Bohlin transformation'' for what I now   call
 the   ``Maclaurin   transformation'', the transformation which yields what I'm calling
 Maclaurin    duality. 
  Arnol'd   used  
``Bohlin transformation''  in his book    \cite{Arnold_Huyg},  which is where  I learned  about this
remarkable duality. 
I sent this early draft to Alain Albouy who   corrected my historical  misunderstandings.
I   copy some of the  history which   Albouy and Zhang \cite{Albouy}  have  unearthed,
restating it  in  an  abbreviated  form. 

\vskip .3cm 
  
1742. Colin Maclaurin   \cite{60}, in section 451 of his Treatise of fluxions,  
   states the duality of the previous section, proving it   in section 875.

1889. Goursat \cite{30}  
rediscovers Maclaurin  duality using  the Hamilton-Jacobi equation and conformal maps.

1889. Darboux \cite{21} extends 
   Goursat \cite{30}  to mechanics in  curved,   multidimensional spaces.

1894. Painlev\'e \cite{70}, inspired by Darboux, asks in   \cite{70}.  He  ``what are all  the transformations which,
together with time reparameterizations,  send one   system onto another?'' He analyzes  previous works on page 17
of this publication.

1896. Painlev\'e \cite{72} answers  his    question  for 
 two degrees of freedom systems.  Darboux's examples, and so Maclaurin  duality,   are one of  Painlev\'e's main cases.   

1900.  Ricci and Levi-Civita \cite{75} develop   covariant differentiation and  the  Ricci calculus,
(also known as the debauch of indices).
In chapter 5, section 4 they cite   Painlev\'e's question 
  as one of their (many) motivations. They do not recall  Maclaurin  duality.   
 
1904. Levi-Civita \cite{52} uses squaring (transformation (\ref{folding map})), a special case of the
Maclaurin   transformation,   to  regularize the planar  restricted 3-body problem.
 
1911. Bohlin \cite{14}
re-presents and perhaps rediscovers  the squaring  transformation of Levi-Civita and
predecessors   as part of a new method of integration of the Kepler problem. He does not give references.
 
1941. Wintner \cite{Wintner}, a classic book in celestial mechanics, on  page 423,  cites Bohlin (op. cit.)  for his ``elegant method of integration''.   
 
1953. Faure \cite{26,27}   remarks that  squaring 
pulls back  Schr\"odinger's equation for  the planar hydrogen atom to that for the planar harmonic oscillator. More generally, he describes Goursat's results in the   planar quantum mechanical context, including   Maclaurin   duality. He does not provide  any references.

1981.   McGehee \cite{McGehee}, in the process of studying the possibility of regularizing of collisions,    rediscovers
the Maclaurin transform $Q \mapsto Q^{\beta}$ for the case $\beta$  an integer  so that the transform is single-valued.
(See subsection \ref{subsec: metric quotient} below.) 
He references Levi-Civita's work on the Kepler case   $\beta = 2$.  
 
1989. Arnol'd and Vasiliev \cite{Arnold_V} 
 describe Maclaurin duality,  refering to it as Bohlin duality and citing Bohlin (op. cit.)  and Faure (op. cit).
 
 \vskip .3cm 
 
This history  provides  a perfect  example of  a favorite saying of  Arnold: 
  A mathematical discovery   is named after a    person
   only if that  person was not the first to disover it.

 \section{Cones at zero energy}
 \label{sec: cones} 
 
 \subsection{Two-dimensional cones} 
 The {\it cone}  of  angle $2 \pi c$   is the  
metric  space defined by the  Riemannian  metric 
 \beq 
 ds ^2 = d \rho ^2 + c^2 \rho^2 d \theta ^2
 \Leq{cone} 
 on the plane where    $(\rho, \theta) $ are  standard  polar coordinates in the plane. 
The parameter  $c$ occurring in the metric is any fixed positive  real number.  
 The origin $\rho = 0$ represents  the cone point and $\rho$ measures distance from it.
 The length of a curve is the integral of $ds$ over that curve
 so that the length of the  circle  of radius $\rho =1$  centered at the cone point
 is $\int ds = 2 \pi c$ rather than the traditional $2 \pi$. 
 The cone  metric is that of the Euclidean
 plane when $c =1$. Otherwise the \Ri metric becomes  singular at the cone point. 
Another term for the cone of angle $2 \pi c$ is 
 the cone over a circle of radius $c$.  Here is a restatement of part B of theorem \ref{thm 1}.

\begin{theorem}
 The JM metric at zero energy  for the   potential $V = - \frac{\mu}{r^{\alpha}}$, $\alpha \ne 2,  \mu > 0$
 is that of a  cone of angle  $2 \pi c$ where 
 $  c = | 1- \frac{\alpha}{2}|$.
 \label{thm: alpha cone}
 \end{theorem}  
 
 {\sc Proof of  theorem \ref{thm: alpha cone}.}
 
 In polar coordinates $|dq|^2 = dr^2 + r^2 d \theta^2$
 so that, up to a constant,  our JM metric is $ds^2 _{JM} = r^{-\alpha} dr^2 + r^{2 -\alpha} d \theta^2$.
 We look for a change of variables $(r, \theta) \mapsto (\rho, \theta)$
 that puts $ds^2 _{JM}$  it into the standard conical form  (\ref{cone}).   This suggests solving
 $d \rho^2 = \frac{d r^2}{r^{\alpha}}$ or $d \rho = \frac{d r} {r^{\alpha/2}}$. 
Solve  by guessing $r = \rho^{\beta}$ and deriving an equation for $\beta$. We get 
  $dr = \beta \rho^{\beta-1}$, $dr^2 = \beta^2 \rho^{2(\beta-1)} d \rho^2$ so that
  $r^{-\alpha} dr^2 = \beta^2 \rho^{2 \beta -2 - \beta \alpha} d \rho^2$.
 We need this  last expression to be a constant times $d \rho^2$ which requires
 $2 \beta - 2 - \beta \alpha =0$. Solve to find  that  $\beta = 2/(2 - \alpha)$ which we recognize
 from earlier.  The term $r^{2-\alpha}$ in front of $d \theta ^2$ 
 in the metric is equal to $\rho^{\beta(2 - \alpha)} = \rho^2$.
 We have converted  the metric into the form $\beta^2 d \rho^2 + \rho ^2 d \theta ^2 =
 \beta^2 (d \rho ^2 +   \frac{1}{\beta^2}  \rho^2 d \theta ^2)$.  
 The dilation (or substitution) $\rho \mapsto \beta \rho$
 converts this last metric to    $d \rho^2 + \frac{1}{\beta^2}  \rho^2 d \theta ^2$
 which is in our   standard conical form (\ref{cone}) with cone  angle constant $c$  determined to be    $c^2 = 1/\beta^2$.
 QED
 

\subsection{Paper folding and the Maclaurin   transformation}
 
 The change of variables
\beq 
 \psi = c \theta
 \Leq{unfolding}
  converts the conical metric (\ref{cone}) to  $d \rho^2 + \rho^2 d \psi^2$
 which is  the flat metric $du^2 + dv^2$ in the $(u, v)$
 plane when we convert to Cartesian coordinates
  in the standard way:  $u = \rho \cos (\psi),  v = \rho \sin(\psi)$. 
 {\it The $(u,v)$ plane  is  the folding plane.}    We   set 
 \beq
 Q = u + iv = \rho e^{i \psi}
 \Leq{folding variable}
to be our standard complex  variable coordinatizing the folding plane.
 
 The change of variables $(\rho, \theta) \to (r, \theta)$ used in the proof of theorem \ref{thm: alpha cone}
 can be viewed as an intermediate step in the Maclaurin   transformation 
 $Mac(Q) =Q^{\beta}$.
 Writing $q = r e^{i \theta}$  we see that 
 the Maclaurin   transform  $Mac$ in polar coordinates is  $(\rho, \psi) \mapsto (\rho^{\beta}, \beta \psi) = (r, \theta)$.
 We   factor this map as the composition of two  maps, one acting
   on   angles, the other   on the distance coordinates.  Thus  
\beq
 Mac  = F_{JM} \circ F_{fold},  \text{ where }  F_{JM} (\rho, \theta) = (\rho^{\beta}, \theta),  F_{fold} (\rho, \psi) = (\rho,   \beta \psi) 
 \Leq{folding}
 The second map which 
 we   call $F_{JM}$ for ``Jacobi-Maupertuis'' 
 was used in the proof of theorem \ref{thm: alpha cone} to put the  JM metric into the standard conical form.  The 
 first map which we   denote by $F_{fold}$ for ``fold''   
 is the folding map, telling us how to construct the cone out of paper.

 As $\theta$ varied over $[0, 2 \pi]$
 our  scaled angle $\psi$ varies over the interval $[0, 2 \pi c]$.
 If $c < 1$  this interval of angles
 is less than $2 \pi$.  This gives us a  paper folding interpretation
 of our  cone. Remove the sector $2 \pi c < \psi < 2 \pi$,
 leaving the sector   $0 \le \psi \le 2 \pi c$ 
 bounded by the  two rays $\psi = 0$ and $\psi = 2 \pi c$.  Since   $\theta$ is defined modulo $2 \pi$
  we must  think of $\psi$ as defined modulo $2 \pi c$  to recover the cone.
  This requires gluing  the bounding rays to each other  by identifying points   $(\rho, 0)$ to points $(\rho, 2 \pi c)$.  
  We have made our cone with paper, scissors and glue. 
  When $c = 1/2$ this is the construction we gave of the Kepler cone.
  
  If $c > 1$ we still have a paper folding interpretation of the cone.  If $1 < c < 2$, slice the plane along some ray,
  and open  the slice so we have now two bounding rays.   Take a new plane and cut out of it a   sector of opening angle $2 \pi (c-1)$ and
  attach its two bounding rays to the sliced plane, gluing one bounding ray to one edge of the slice
  and another ray to the other. If $c = k$ is an integer greater than $1$   we must 
    glue  in series $k$ slit planes, resulting in an
  object   modelling  of a k-fold  branched cover, making a total angle of $2 \pi k$ about the
  cone point and then we glue in the final half of the last slit to the
  first part of the initial slit plane.  If $c$ is not an integer and $k < c < k+1$ we glue in a 
  final sector of opening angle $2 \pi (c-k)$ before we close the object up.
  
  The folding map, like the Maclaurin   transformation,  is multi-valued. To make it into an honest map
  we can  interpret $\psi$ as a real variable, not an angle.  Topologically this corresponds
  to understanding  that $(0, \infty) \times \R$
  is the  universal cover of the punctured plane, or cone minus the cone point,
  and using  $(\rho, \psi) \in (0, \infty) \times \R$ to coordinatize this universal cover.  
  The fundamental group of the punctured plane is $\Z$,
  (the fundamental group of the circle) and it acts on the universal cover so that   $k \in \Z$
  acts by $(\rho, \psi) \mapsto   (\rho, \psi + k (2 \pi c))$.   
Our paper work, the   process of slits and gluing,  is a paper-and-scissors way of realizing
  that the cone is the metric quotient of this universal cover by this action of the fundamental group.
 (Endow the  universal cover  with the conical metric $d \rho^2 + \rho^2 d \psi^2$  , remembering to
 take  $\psi$  as a real not an angular variable.)

\subsection{Cylinders:  the   case   $\alpha = 2$}.
\label{subsec: cylinder}

The  cone construction  breaks down when the radius $c$
of the circle is zero which occurs when $\alpha =2$. 
We can rescue geometry    by  carrying out the JM metric analysis
for the corresponding potential $V = - \mu/ r^2$.    The  JM metric is
   $\frac{2}{r^2} (dr^2 + r^2 d \theta ^2) = 2 ( (\frac{ dr}{r} )^2 + d \theta ^2)$.
Since  $\frac{d r^2} {r^2} = d Log (r)^2$ we set $\rho = Log(r)$.  In $(\rho, \theta)$ 
variables the JM metric is (up the constant scale $2$)
 $d \rho ^2 + d \theta ^2$.  This is the  induced metric on the cylinder $x^2 + y^2 = 1$ in Euclidean
 $\R^3$ endowed with standard $x, y, z$ coordinates, upon setting  $\rho =z$.  The metric is flat, rotationally symmetric
 and admits the additional isometries $(\rho, \theta) \mapsto (\rho+ \rho_0, \theta)$
 and $(\rho, \theta) \mapsto (- \rho, \theta)$.   This translational isometry
 is a kind of limiting ghost of the dilations for the cone.  
 
 The cylinder  has three types of geodesics: circles $\rho = $ const,
 generators $\theta = $ const, and then the typical geodesic, 
 a helix.  All of these types can be achieved by rolling the cylinder along  a plane
 after having painted a line on the plane.    Back on the Newtonian plane
 the circular geodesics are circles, the generators become rays
 and the  helices become logarithmic spirals $q = e^{k t} e^{i \omega t}$,
 $k, \omega \ne 0$ real parameters.

\section{Geodesics on  cones.}

Cones have   two types of geodesics:  those that hit the cone point,
and those that don't.    The geodesics which hit the cone point are the rays, $\theta = const$, also known as the generators.  What does  a non-collision geodesic look like?   
It comes in from infinity  asymptotic
to one generator $\theta = \theta_-$, gets close to the cone point, then returns to infinity asymptotic to another  generator
$\theta = \theta_+$.  The angle between these two generators is
\beq
|  \theta_+ - \theta_- | = \frac{\pi}{c}
\Leq{scattering angle} 
and is  the same for all non-collision geodesics.
We call the quantity (\ref{scattering angle})  the ``scattering angle'' of the 
geodesic, this,  despite the fact that  it is really a real number,
not an angle.  For example,  if $c =1/4$
the scattering angle of non-collision geodesics   is $4 \pi$.
This `$4 \pi$' gives us two pieces of information,
First, the fact that $4 \pi \equiv 0 $ mod $2 \pi$ signals 
that the incoming asymptotic ray and outgoing asymptotic
ray are the same.   Second,  the fact that  $4 \pi = 2 * 2 \pi$
  tells us that a  non-collision
 geodesic on this cone   will cross 
every ray exactly twice, except for its incoming-outgoing asymptotic ray which 
 it hits exactly one.  

Allow us to be somewhat formal regarding what we've just said.  
\begin{definition}  
The scattering angle  of a collision-free curve on the
cone   is the measure of the set of rays, counted
with mulitiplicity, crossed  by that curve.  
\end{definition} 

\begin{lemma} 
The  scattering angle   of a
  non-collision geodesic on the cone 
of angle $2 \pi c$  
is 
  given by   equation (\ref{scattering angle}).  
\label{lem: angular measure}
\end{lemma}

{\sc Proof.} Represent the geodesic  by a line $\ell$ in the folding plane,
with  polar
coordinates are $\rho, \psi$.  Being a straight line, its scattering angle as a curve on
the  
folding plane   is $\Delta \psi = \pi$, representing the fact it that hits
`half' the rays in the folding plane, namely those lying in the half-plane containing
$\ell$ whose bounding line is parallel to
the line.   Alternatively,  parameterize $\ell$ by arclength $s$ and   express it in
polar coordinates as a function of arclength:  $( \rho(s), \psi(s))$.  
Set $\psi_{\pm} = \lim_{s \to \pm \infty} \psi(s)$.  Then   
$| \psi_+ - \psi_-| = \pi$ as we can see by  drawing a picture .  
The angular function $\psi (s)$ is strictly monotonic with respect to $s$.  
Now let $(\rho, \theta)$ be polar coordinates on the cone.  From  the relation $\theta = \psi/c$.
we see that our geodesic is given by $(\rho(s), \frac{1}{c} \psi(s))$
in the  cone coordinates.   The angle $\theta$ is also strictly monotone. 
It follows that  the scattering angle of the geodesic is 
$\Delta \theta = \pi/c$.

\begin{figure}
\scalebox{0.4}
{\includegraphics{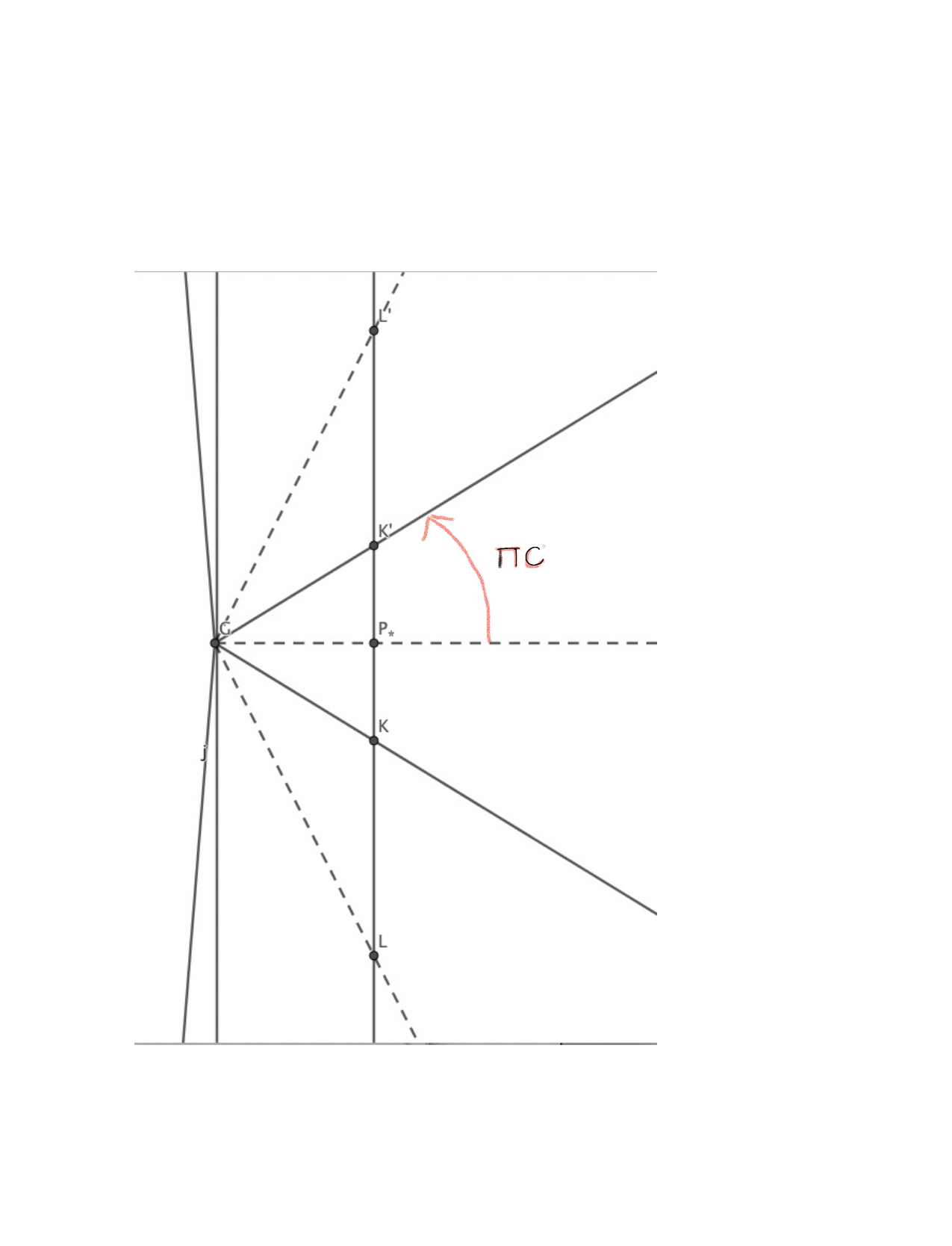}}
\caption{The fundamental sector can be used to tile over the half-plane containing the model line.
$P_*$ represents the closest point on  on the corresponding geodesic to the cone point. Four of the five intersections of  the line with the
rays forming  the boundaries and bisectors of the tiles are marked, L, L', K, K'.  These
correspond in pairs to self-intersections of the corresponding geodesic on the cone over the
circle of radius c. There are $k+1 =3$ tiles and correspondingly $k =2$ self-intersections.}  
\label{fig: conetilingC}
\end{figure}

\begin{exercise}  Show that if the scattering angle is 
greater than $2 \pi$ then the non-collision geodesic   
self-intersects.   Let $\theta = \theta_*$ be the ray containing the point $P_*$
 on the geodesic closest to the cone point.  Show that the self-intersections
  are alternately located along the ray $\theta = \theta_*$  and its antipodal ray $\theta_* + \pi$.
Let $k$ be the integer part $1/2c$ so that  
the scattering angle is  $k 2 \pi + 2 \pi r$ with  $0 \le r < 1$.
Show that if the fractional part  $2 \pi r$ of the scattering angle is nonzero then
the geodesic has exactly $k$ self intersections.   See figure \ref{fig: conetilingC}.
\end{exercise}

\subsubsection{Zero energy scattering for power laws.}  
Lemma \ref{lem: angular measure}  
combined with the Maclaurin transformation 
gives us the same scattering angle  for the zero-energy solutions to
our central force law.  We continue to be formal around definitions.  
\begin{definition}  
The scattering angle  of a collision-free curve on  the Newtonian plane  is the measure of the set of rays, counted
with mulitiplicity, crossed  by that curve.  
\end{definition} 
\begin{lemma}
\label{scattering zero energy} 
 The scattering angle  of any zero-energy solution
 the power law central force problem with exponent $\alpha$
 is given by equation (\ref{scattering angle})
 where $c$ and $\alpha$ are   as before:   $c = 1- \alpha/2$.  
\label{lem: scattering measure 2}
\end{lemma} 

{\sc Proof}  Recall the Jacobi-Maupertuis map $F_{JM}$,
  $(\rho, \theta) \mapsto (\rho^{\beta}, \theta) = (r, \theta)$
in the factorization (\ref{folding}) the Maclaurin transformation.
This map $F_{JM}$ takes  non-collision geodesics
on the cone of angle $2 \pi c$ to non-collision zero energy solutions to Newton's equations
with exponent $\alpha$.  It  leaves  
$\theta = \psi/c$ unchanged and  acts on the distances in a nice
monotonic way.     Consequently the count for
the number of rays hit by the solutions to Newton's equations
is identical to the number of rays hit by its corresponding geodesic.
 QED

  \subsection{Integer powers and Metric Quotients}   
 \label{subsec: metric quotient}
 
 The scattering angle formula (\ref{scattering angle})
shows us  that the incoming and outgoing asymptotes of non-collision geodesics
 are the same precisely  when $\frac{1}{c} =n$ for some  positive {\it even}  integer $n$.
Recall that  $\frac{1}{c} =  \beta$ where $\beta$ is
 the exponent  of the Maclaurin transformation
 $Q \to Q^{\beta}$.      
 If $\beta = n$ is a positive integer (even or odd)   then, and only then, does
 the Maclaurin   transformation $Q \mapsto Q^n$
become a  single-valued map from the folding plane to the Newtonian plane.
 Under this map, two points $Q, Q'$ get mapped to the same point $q$ under this map if and only if they
 differ by scalar multiplication by some n-th root of unity,
 which is to say that  points on the fibers of the MacLaurin transform
   are related to each other by   integer multiples of 
 rotation by $2 \pi /n$.  It follows that  the map $Q \mapsto Q^n$ defines an isomorphism
 between  the quotient space  $\C/ \Z_n$ and the Newtonian plane $\C$.
 If we take this quotient $\C/ \Z_n$ to be the metric quotient then we have an algebraic
 description of our  cones with opening angle $ 2\pi c,  c = 1/n$, $n =1,2,3, \ldots$.   
  
Some words are in order regarding the adjective ``metric'' in the phrase  `metric quotient space'
used just above and used earlier   to describe the Kepler cone as $\C/\pm 1 = \C/\Z_2$. 
If  a finite group $G$ acts isometrically on a metric space $M$, then the  quotient space
$M/G$ naturally inherits a metric by declaring the distance between two points
in the quotient to be the distance between the corresponding orbits in $M$.  
To obtain the Kepler cone as a metric quotient
we take  $G = \Z_2 = \{1, -1\}$ acting on  $M = \C=  \R^2$ by scalar multiplication,
putting the  standard Euclidean metric on $M$.   To obtain these other cones $\C/\Z_n$
we use the action by rotation by angles $2 \pi j / n$, $j =0, 1, 2, \ldots$ as realized by
scalar multiplication by the group of nth roots of unity.  

 On any   metric space $(M,d)$ we have the  notion of a 
curve being a  ``geodesic''.  To define geodesics on $M$  first   assign a length (possibly infinite!) to any continuous
curve $c:[a, b] \to M$  in a metric space by using a process of successive `polygonal' approximations.
\footnote{The approximating lengths are then the  sums $\Sigma d(c(t_i), c(t_{i+1}))$ where $a = t_0 < t_1 < \ldots t_n =b$
is a partition of the interval $[a,b]$. The length itself is the lim sup of these approximating lengths.}   
Then   declare that a continuous  curve is a minimizing geodesic between two points if its length is
equal to  the distance (i.e $d(c(a), c(b))$)
between its endpoints.  And declare  a curve to be  a (not-necessarily minimizing) 
geodesic if all of its sufficiently short subarcs
are minimizing geodesics (between the endpoints of these subarcs).   See Burago et al \cite{Burago}
for details. 

The metric quotient operation
$M \to M/G$ maps geodesics to geodesics.  The geodesics on $\R^2$ are the straight lines.   So the geodesics on these `integer cones'' $\C/\Z_n$ are the quotient images of the straight lines in the $(u,v)$ plane. 
The metric quotient
$\C/\Z_n$ is isometric to a cone  of angle $2 \pi/n$ since a  fundamental domain for the
action of $\Z_n$ is the  sector $0 \le \psi \le 2 \pi/n$.

\section{Scattering and starbursts:  Implications for Non-zero energies}

\begin{figure}[ht]
 \scalebox{0.3}{\includegraphics{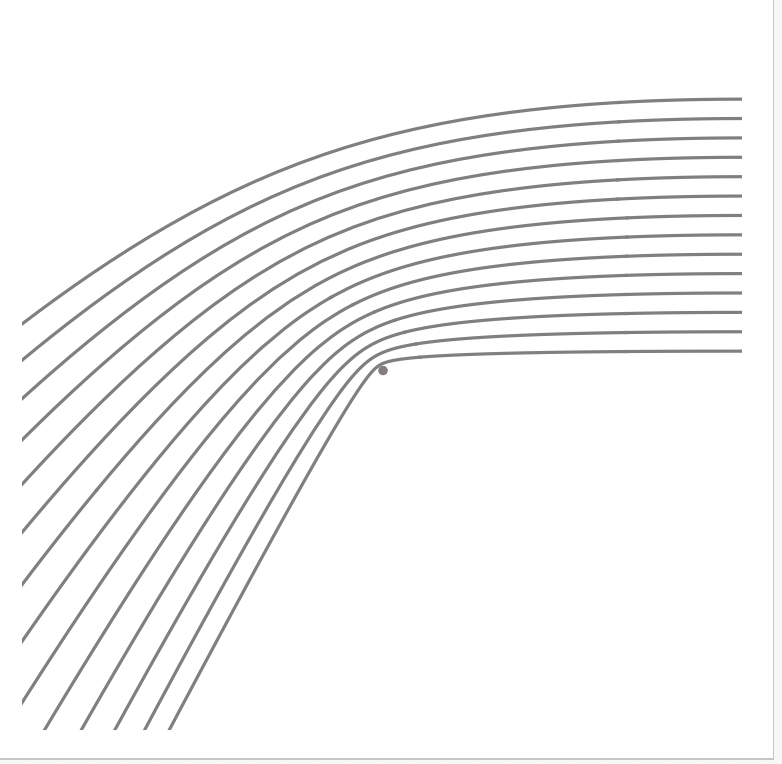}}
\caption{Positive energies.  Scattering a half beam  with power law $r^{-1/2}$ so  $\pi/c = \frac{4 \pi}{3}$.}
\label{fig: twothirds_scattered}
\end{figure} 
\begin{figure}[h]
\scalebox{0.3}{\includegraphics{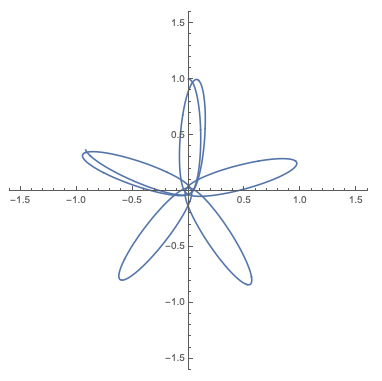}}
\scalebox{0.3}{\includegraphics{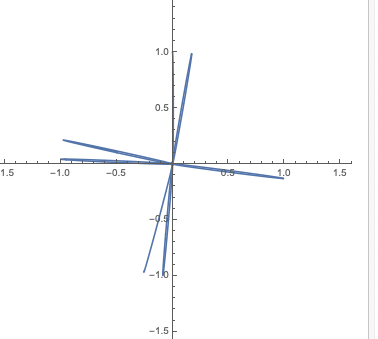}}
\caption{Negative energy solutions, or `rosettes' for $\alpha = 2/3$ in
which case $\pi/c = \frac{3 \pi}{2}$. The rosettes  converge to a $\pi/c$ starburst as the angular momentum tends to zero.}  \label{starburst}
\end{figure}

Zero-energy solutions with their scattering angle (\ref{scattering angle})
continue to have relevance at non-zero energies.  
 Figures \ref{fig: twothirds_scattered} and \ref{starburst}
summarize how.   They indicate families of such solutions converging
to a curve suffering a  collision.  In both cases the  limiting curve exhibits a change of direction
at collision, this change equaling the  scattering angle (\ref{scattering angle}).  These two  figures illustrate
 \begin{theorem}  Let $\gamma_{\eps} (t)$ be a family of non-collision solutions to the $\alpha$-power law problem
(equation (\ref{eq: C force 1}))having   energy $E$
and suppose that  $\gamma_{\eps}$ converges as $\eps \to 0$  to a curve $\gamma_*$
having a collision.    If   $0< \alpha < 2$  then  $\gamma_*$   parameterizes a piecewise linear curve
whose only vertices  (points with changing directions of travel) are at the collision (origin)  and   the Hill boundary (if
there is one).  Away from the origin the limiting curve satisfies Newton's equations. 
Upon each  collision the direction of $\gamma_*$ changes, turning
by an  angle equal to the scattering
angle $\pi/c$ with $c = 1 - \frac{\alpha}{2}$ as usual.
\label{thm: limit}
\end{theorem} 

When $E \ge 0$  the limiting curve sweeps out  the union of two rays which make  an angle equal
to the scattering angle.  
When $E< 0$ the solutions are bounded, lying in 
 the Hill disc $r \le (\mu/E)^{1/\alpha}$ and then  the  limiting curve 
 is a `starburst' concatenation of brake orbits, with successive orbits
 making an angle equal to the scattering angle with each other.
  Figure \ref{fig: twothirds_scattered}   depicts the positive energy case
  by exhibiting  an incoming half-beam of particles 
 under the influence of an $\alpha = 1/2$
 power law potential. In the limit as  the  solutions in the  beam tend to the 
 origin
 their  scattering angle limits to  $\pi/c = \frac{4}{3} \pi$ as described in the theorem.
 Figure \ref{starburst} depicts the negative energy case
 by showing  two negative energy
 solutions  to an $\alpha = 2/3$ power law problem.  The second solution
 comes very close to the origin and successive near-brake solutions almost
 make the angle $\pi/c = \frac{3}{2} \pi$ with each other.   
The `starburst' aspect of the theorem seems to be a   previously unobserved fact regarding power law central force laws. 
  
 {\sc Setting up a proof.}  Power law problems,  like all  central force problems,  admit
  rotational symmetry and conservation of angular momentum 
\begin{eqnarray}
J &=& q \wedge v \\
&=& r^2 \dot \theta.
\label{J1_2}
\end{eqnarray}
 Rotational symmetry means that
if $q(t)$ is a solution then so is any rotate  $e^{i \theta_0} q(t)$ of this
solution.  The first expression for the angular momentum $J$ in equations (\ref{J1_2}) uses the two-dimensional version of the  cross product $(q, v) \mapsto q \wedge v = (x,y) \wedge (v_x, v_y) = x v_y - y v_x$.  
The second expression  for $J$ in equations (\ref{J1_2})uses  polar coordinates  $(r, \theta)$. 
Conservation of angular momentum means that $J$ remains  constant along solutions. 
Solutions to central force problems  travel along rays through the origin if and only if $J =0$
as one can see from the second expresson for $J$.  

\begin{figure}[ht]
\scalebox{0.4}{\includegraphics{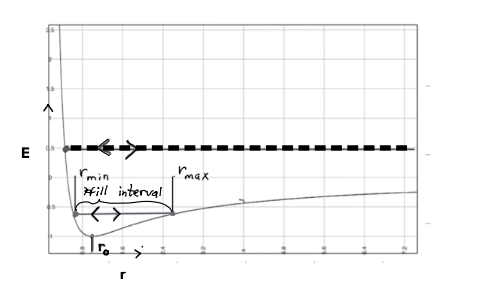}}
\caption{The graph of  the effective potential $V_J (r) = J^2 /r^2 - \mu/r$ for the Kepler problem.  
I took $J =1, \mu =2$. 
The horizontal lines represent fixed energy levels.  The thin solid one is negative and $r_0$
is the circular orbit for this angular momentum.  The dotted thicker line  represents a positive energy value with corresponding unbounded, or scattering orbit.  }
\label{fig: KeplerEffective}
\end{figure}

If a particular solution  to a  power law problem 
with $\alpha \in (0,2)$  tends  to the  the origin with time
then the solution must have   $J = 0$.
We can show this by re-expressing Newton's equations  polar coordinates.
The kinetic energy is given by  
$K = \frac{1}{2} (\dot r^2 + r^2 \dot \theta ^2)$ from which it follows that
$$E = \frac{1}{2} \dot r^2 + \frac{J^2}{2 r^2}  - \frac{\mu}{r^{\alpha}}$$
 since $r^2 \dot \theta ^2 = \frac{J^2}{r^2}$.   The  sum of the last two
 terms in this rendiction of energy is called the   ``effective potential'' $V_{eff}$: 
 $$V_{eff} (r; J) = \frac{J^2}{2 r^2}  - \frac{\mu}{r^{\alpha}}$$
 (Newton's  equations, rewritten in polar coordinates, are   $\ddot r = -\frac{d}{dr} V_{eff} (r ; J)$
together with  $\dot \theta = \frac{J}{r^2}$ and $\dot J = 0$.) 
 Since $\frac{1}{2} \dot r^2 \ge 0$ we have that 
 $V_{eff} (r; J ) \le E$.   It follows that for fixed $E, J$
 with $J \ne 0$  the radial variable $r$   lies in the `Hill interval''  $\{r:  V_{eff} (r; J) \le E \}$.
 The salient fact for us is that when $J \ne 0$ and  $\alpha$ is  in our range 
 the left endpoint of the  Hill interval is a positive number $r_{min} (E, J) >0$
 and that 
 \beq 
 r_{min} (E, J) \to 0 \iff J \to 0
 \Leq{eq: rmin}
   for $E$ fixed or lying in a bounded interval. 
 (If $E \ge 0$ then   the Hill set   is of the form $[r_{min} (E, J), \infty)$
while  if $E <  0$ it has the form $[r_{min}(E, J), r_{max} (E, J)]$
where $r_{max}(E, J)$ is finite.)  We can  see this  fact regarding the Hill interval by graphing  
 $V_{eff} (r; J)$ as a function of $r$  for the parameter  $J \ne 0$.  See  figure \ref{fig: KeplerEffective}.  
 We have, for  $\alpha$ in our range $(0,2)$,   that
 $\lim_{r \to 0} V_{eff} (r;J) = + \infty,   \lim_{r \to \infty} V_{eff} (r;J) = 0$
 and that $V_{eff}( r; J)$ has only one critical point $r = r_* (J)$ in the range $0 < r < \infty$
 and this critical point is a   global minimum 
 whose value is negative.     The structure of the Hill interval follows.  
 
 \begin{exercise}  For $\alpha = 2$ show that there are solutions
 with $J \ne 0$ which are  logarithmic spirals, limiting to the origin. 
 Show that whenever $\alpha \ge 2$ there are solutions to the
 power law problem having $J \ne 0$ but which   limit to the origin,
 spiralling in as they go.  \end{exercise}

 {\sc A circular solution space. }    Fixing  $E$ and $J$  
defines  a  family  ${\mathcal F} (E, J)$  of solutions  mapped to itself  by rotations.
For power law potentials in the range $0 < \alpha < 2$ the family
${\mathcal F} (E, J)$   consists of a single solution modulo rotations
and time translations, provided $J \ne 0$ and provided the family is not empty.
Theorem \ref{thm: limit} concerns 
  $\lim_{J \to 0}  {\mathcal F} (E, J)$. 

\vskip .3cm

    {\sc Proof of the theorem.} 
 
Write $r_{min} (\eps) = min_t | \gamma_{\eps} (t)|$
  for the distance from the origin  of the solution $\gamma_{\eps}$.
By hypothesis $r_{min} (\eps) \to 0$
  with $\eps$   so that  according to equation
  (\ref{eq: rmin})  we have that its angular momentum $J \to 0$ with $\eps$.  
   It follows that away from collisions our limiting curve $\gamma_*$  
obeys   Newton's equations and has $J = 0$ and so
lies along rays.  That ray may  change only upon  collision.
 Proving the theorem
amounts to showing that the ray does change at each collision and that
the change in angle is   precisely the scattering angle.

To verify the change in angle we use an additional    space-time 
 scaling symmetry 
  \begin{equation}
 \label{scaling}
\delta_{\lambda}:    q  (t) \mapsto  \lambda q (\lambda ^{-\nu} t);  \nu = \frac{\alpha}{2} + 1.
 \end{equation} 
 enjoyed by $\alpha$- power law  problems.     
This transformation  takes solutions to solutions,
a fact which follows from the homogeneity of the potential $V = -\mu/r^{\alpha}$.
This  scaling  transforms velocities
by $v \mapsto \lambda^{-\alpha/2} v $ and energy and angular momentum
by $E \mapsto \lambda^{-\alpha} E,  J \mapsto \lambda ^{1- \alpha/2} J = \lambda^c J$.  
 Thus the scaling 
 defines a scaling  isomorphism   $\delta_{\lambda}:  {\mathcal F} (E, J) \to {\mathcal F} (\lambda^{-\alpha} E,  \lambda^c J)$ between solution spaces.
 Here we have used the notation $ {\mathcal F} (E, J)$ introduced in   the remark
 immediately above for the space (a circle) of all solutions
 having energy $E$ and angular momentum $J$.  Scale invariant properties of these curves, such as the
 angle between successive  perihelion and   apihelion  when $E < 0$, are
 necessarily preserved by the scaling.    Now, set $J = \eps$, imagining $\eps \to 0$
 and set  $\eps = \lambda^{-c}$ so that $\lambda = 1/\eps^{1/c} \to \infty$ as $\eps \to 0$. 
 The scaling transformation takes  ${\mathcal F} (E, \eps)$ to ${\mathcal F} (\lambda^{-\alpha} E,  1 ) =
 {\mathcal F} (\eps^{\alpha/c}  E,  1 )$.  It follows that the scattering angle associated to
 fixing $E$ and letting $J \to 0$ is the same as the scattering angle associated
 to fixing $J= 1 \ne 0$ and letting $E  \to 0$.      But this is the scattering angle for 
 ${\mathcal F} (0, 1)$  as described in lemma \ref{lem: scattering measure 2}.
 
 QED
      
  \subsection{ Scattering off  power law potentials}
   
\begin{figure}[ht]
\scalebox{0.4}{\includegraphics{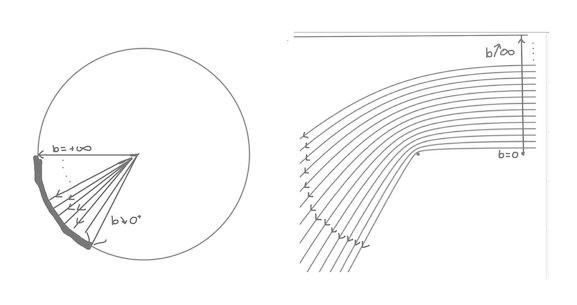}}
\caption{The right panel shows a half beam being scattered by an attractive central force  power law $f(r) = r^{-1/2}$.
The left panel shows the arc of the circle resulting from  the scattering map applied to this half beam, 
its image being half of the scattered  image of the full beam. }  
\label{fig: halfbeam}
\end{figure}

 Figure \ref{fig: halfbeam} is a redrawing of   figure \ref{fig: twothirds_scattered},
 to place it into the   the context of
 classical scattering.  The figure indicates facts I learned   from    \cite{KnaufKrapf} and which I have  summarized below as proposition \ref{prop: scattering}.
Understanding and rederiving these scattering  fact provided the   inspiration
   for writing  this paper.  For further reading on classical scattering
   see  the  chapter on scattering in \cite{Knauf}  and the
   primary source \cite{Rutherford}.

Figure   \ref{fig: halfbeam}  depicts a half-beam scattered by
  an $\alpha = 1/2$ power law central force.   The parameter $b$ of
  the figure is a    transverse parameter   to the  rays comprising the beam
  and $b$  
  labels the   trajectories making up the   half-beam. 
 This parameter  $b$  coordinatizes a line    placed very far from the origin
nearly  orthogonal to the  incoming beam.  As $b \to \infty$ the corresponding
trajectories recede to infinity, barely affected by the central force. 
 As $b \to 0$ the trajectories tend to collision   and so we are in the situation of  lemma \ref{lem: scattering measure 2}.
  The limiting trajectory is  deflected from that of a straight line by
  an angle  $|\pi/c - \pi | = \pi/3 $  as $b \to 0$.  The  $\pi/c = 4 \pi/3$ arises from the scattering angle of lemma \ref{lem: scattering measure 2}
  and   we have to  subtract $\pi$  because, in our definition of `scattering angle'
  an undeflected straight line would have `scattering angle of $\pi$' and we are now
  measuring so that `undeflected' means angle $0$.

  The labeling parameter $b$ is known as the impact parameter.  Its magnitude is the distance of a line
  from the origin if the force had been  turned off so that no deflection occured.  For any $\alpha$-power
  law with $0 < \alpha < 2$ the scattering of a positive energy half-beam shares the  properties just
  described in the paragraph above for $\alpha =1/2$.  
  Write $f(b)$ for the angle of deflection of the trajectory labelled by $b$. 
  We have  $f(b) \to 0$ as $b \to \infty$ and $f(b) = \frac{\pi}{c} - \pi$ as $b \to 0^+$.
  In between $f(b)$ is continuous.  With a bit more work one can show $f$ is monotone $b$
  varies over $(0, \infty)$.
  It follows that 
  \beq f((0, \infty)) = (0, A)
  \text{ where } A =  \frac{\pi}{c} - \pi
  \label{A}
  \eeq
  This arc $(0, A) \subset \bS^1$ describes the  arc of outgoing    rays which arise from the  scattered half-beam.
  
Now consider scattering a whole beam, so that $b$
ranges over the real line $\R$.  
The deflection angle $f(b)$ is continuous in $b$  except at $b =0$.   By reflectional symmetry of  central force dynamics we have that
$f(-b) = -f(b)$.   Trajectories  passing above the center, so with  $b > 0$ in our figure,    swerve to the left
while  trajectories passing  below  the center, so with  $b < 0$,  swerve to the right but the magnitude of 
these two deflections is same for $b$ and $-b$.   It follows that
$f(\R \setminus \{0\}) = (-A, 0) \cup (0, A) \subset \bS^1$
where $A$ is as in equation (\ref{A}) just above.  We have proved: 
\begin{proposition} The image of  a beam upon scattering by an $\alpha$-power law,  $0 < \alpha < 2$,  
is an arc of measure  $2A = 2 \pi \frac{\alpha} {2 - \alpha}$
with its midpoint deleted.  The midpoint corresponds to the  direction  of the incoming beam
and to $lim_{b \to \pm \infty} f(b)$.  The endpoints of the arc correspond to  the deflection
suffered by the two collision  limits $lim_{b \to 0^+} f(b)$
and $lim_{b \to 0^+} f(b)$.
\label{prop: scattering} 
\end{proposition}

If  $\alpha < 1$ then $2A < 2 \pi$ which means
the set of scattered rays does not cover the circle.
The set of missing rays is an arc centered at ``complete backscatter''
meaning  the direction $\pi$ which is  opposite to the incoming beam.  The
arclength of the missing rays is $2 \pi -2 A$.  On the other hand, if $1 < \alpha < 2$
then $2A > 2 \pi$ which means that some outgoing directions are covered more than once:
there are two or more incoming rays which scatter out along a particular outgoing direction.

  \begin{exercise} Show that impact parameter $b$, angular momentum $J$
  and energy $E$ are related by $b = J/\sqrt{2E}$ when we are in the
  scattering regime, $E > 0$.
  \end{exercise}

\subsection{Rutherford's scattering}   
The Coulomb problem is the Kepler problem but with the sign of the coupling constant
reversed so as to make the force repulsive:
\beq
\ddot q =  + \mu q/ r^3,
\Leq{Coulomb} 
with  $\mu > 0$.   It is basic to electromagnetism
where $\mu$ becomes  proportional to the product of the charges involved, the one at
the center and the incoming one.  We have $\mu > 0$ for like charges and $\mu < 0$
for opposite charges.  Like
charges repel and opposite charges attract.  Scattering in the Coulomb problem
is almost identical to scattering in the Kepler problem.  The sign of the deflection
of rays is reversed between the two.  

\begin{figure}[h]
\scalebox{0.3}{\includegraphics{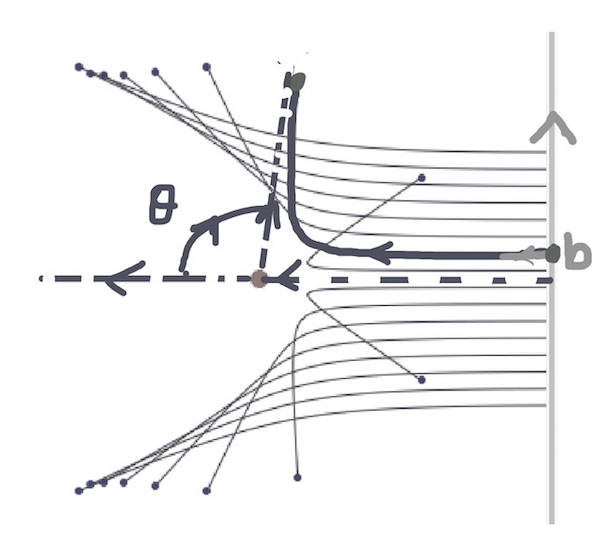}}
\scalebox{0.3}{\includegraphics{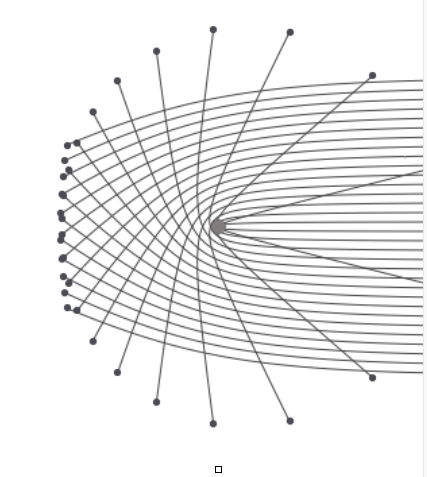}}
 \caption{Rutherford Scattering. The variable $b$ is the impact parameter and $\theta$ coordinatizes the direction 
 of outgoing rays.} 
 \label{fig:Rutherford}
\end{figure}

    The discovery  of the nucleus  relied    
upon  Rutherford's  analysis \cite{Rutherford} of Coulomb scattering.  
   See figure \ref{fig:Rutherford}. 
 His lab assistants Geiger (of the counter) and Marsden had been doing
 experiments in which they directed 
a high-speed beam of alpha particles (Helium nuclei) at gold foil.  They measured the   directions of the outgoing particles scattered by the foil and  were surprised that  a   measurable fraction (1/20,000)  suffered over a ninety degree change in direction.
A few particles suffered 
a nearly   complete recoil.    According to  the dominant  `plum pudding model' of the time,  none should be coming back.  
 To   explain their  experimental  results   Rutherford   replaced the  plum pudding mush  making up the gold atom's center  with  a concentrated point-positive charge containing nearly all the 
 nearly all    the  mass of a gold atom  is concentrated in a   point-like positive  charge  at the  center -the
 nucleus.  A much lighter and  diffuse cloud of   negative charges
(electrons) were to surround the nucleus.  Alpha particles were known to have positive charge, like the nucleus.
The idea was that the alpha particles  penetrated the electron cloud with ease
at which point they had to contend with the strong repulsive  Coulomb force of the nucleus.

 Rutherford computed    
\beq
\theta = f(b) =  2 arctan(\frac{\mu}{2E b})
\label{Rutherford scattering}
\eeq
for  Coulomb-Kepler scattering, that is, for equation (\ref{Coulomb}) regardless of the sign
of $\mu$.  (Note: we just switched sign conventions
compared to our original use of the coupling constant in Kepler.)  
See \cite{Rutherford}, or \cite{Knauf}, particularly   p 286, eq (12.3.3),  or \cite{Landau-Lifshitz}  for this computation.  


 
\begin{figure}
\scalebox{0.1}{\includegraphics{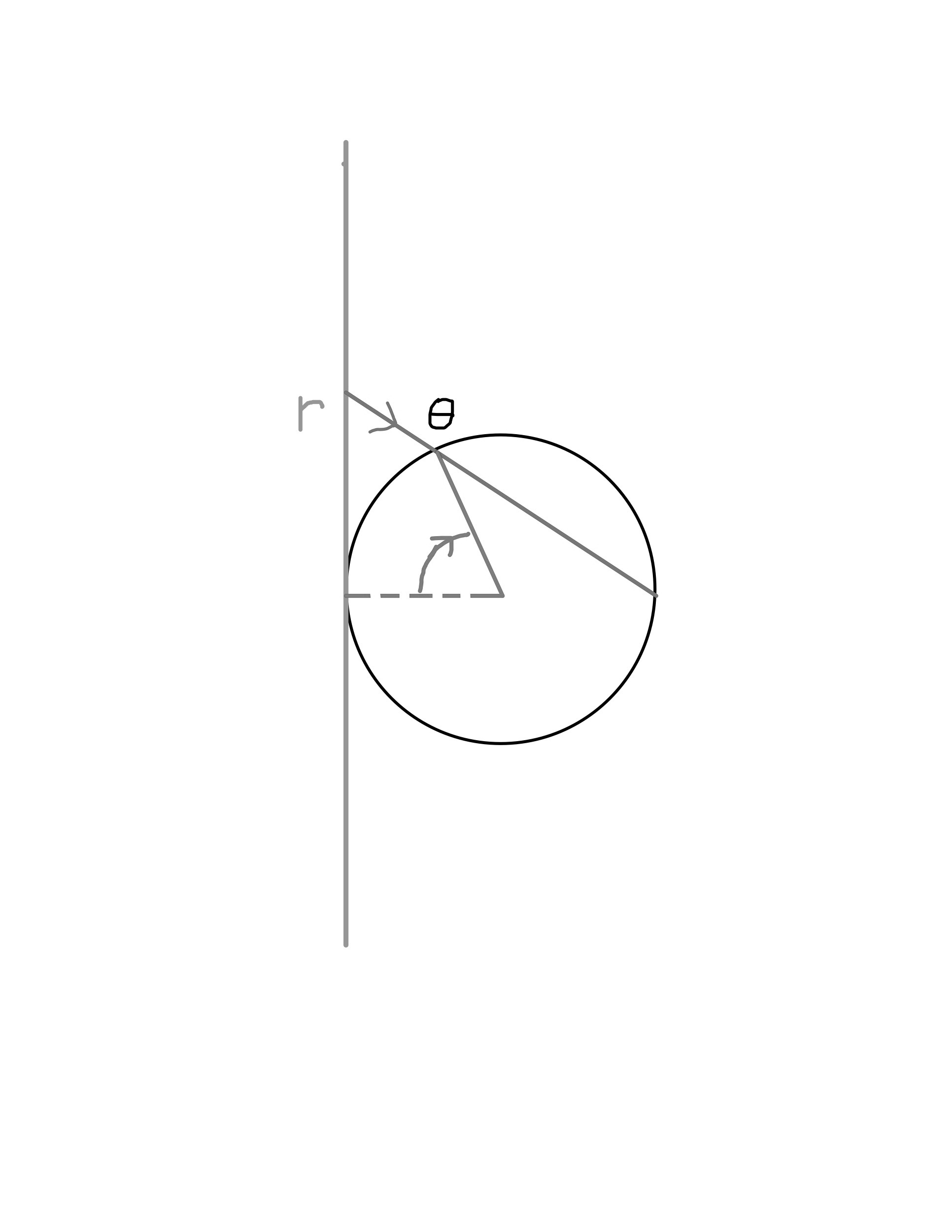}}
 \caption{Stereographic Projection.  Take $r = b$ on the line.} 
 \label{stereoProj}
\end{figure}

It is a beautiful fact that the Coulomb-Kepler scattering map  $f(b)$ 
of formula  (\ref{Rutherford scattering}) is a  stereographic projection $\R \to \bS^1$.   See figure \ref{stereoProj}.
The expression for the stereographic projection $r \to \theta$
indicated in the figure is $\theta =  2 arctan(1/r)$.  
Upon setting  $r = 2Eb/\mu$ we find that Coulomb scattering and stereographic
projection are the same map.  

This relation to stereographic projection makes it clear that the Rutherford scattering map
covers the circle of outgoing directions  {\it exactly once}.  We have   $f(\R \setminus \{0\}) = \bS^1 \setminus \{0, \pi\}$.
We can complete the map   by taking limits:   $f(\infty) =0$ and $f(0) = \pi = - \pi$
The fact that the map extends by continuity for the whole compactified real line
is essentially equivalent to the single-valuedness of Levi-Civita regularization.  
An analogous extension holds whenever $\beta = n$ in the Maclaurin transformation,  as per subsection \ref{subsec: metric quotient}.  If we complete the real line to a circle by identifying
$+ \infty$ to $-\infty$ then  the scattering map 
becomes a degree $n$ analytic map from the circle to itself.   For   details see \cite{KnaufKrapf}.

\section{Deriving the Jacobi-Maupertuis principle.}
\label{sec: JM proofs}  

We will give two derivations of the JM principle,
one  Lagrangian, the other  Hamiltonian.  The first has the
advantage of a direct variational viewpoint but the disadvantage of not easily
yielding the relation between time and JM arclength.

\subsection{Lagrangian Mechanics.}   
The Lagrangian formulation of   classical  mechanics
is  based on choosing  a function $L: \R^n \times \R^n \to \R$ ,  written  
$L(q^i, v^i)$ so that $(q^1, \ldots , q^n,  v^1, \ldots , v^n)$
are used as coordinates of  $\R^n \times \R^n$.
The Euler-Lagrange [EL] equations for $L$ are the ODEs
$$\frac{d}{dt} (\dd{L}{v^i}) = \dd{L}{q^i}$$
supplemented by
$$\dot q^i = v^i.$$ To recover Newton's equations as EL equations take  
$$L (q, v) = \text{ kinetic} - \text{potential} =   K(v)  -  V(q)$$

\begin{exercise}  Show that our  Newton's equations
(\ref{Newtons}) are the   Euler-Lagrange equations
for  $L(q, v) = \frac{1}{2} \| v \|^2 - V(q)$
where $q = (x,y),  v = (v_x, v_y)$ and $(q^1, q^2, v^1, v^2) = (x, y, v_x, v_y)$. 
\end{exercise} 

The Lagrangian connects the calculus of variations to mechanics.
Define the action of an absolutely continuous path to
$q: [0, T] \to \R^n$ to be $\int_0 ^T L (q(t), \dot q(t)) dt$.
Consider the ``standard  problem of the calculus of variations'':
to minimize   the
action among all paths joining two fixed endpoints  $q_0, q_1$ in a time
$T$.   The minimizer, {\it if it exists, and if it is twice-differentiable}, must  
satisfy the EL equations and hence Newton's equations.      

\begin{example}  [Riemannian geometry]  A Riemannian metric is given in local coordinates  by
a Lagrangian which is quadratic positive-definite in velocities:
\beq
L(q, v) = \frac{1}{2} \Sigma g_{ij} (q) v^i v^j
\eeq
where the $q$-dependent matrix $g_{ij} (q)$ is symmetric positive-definite.
$2L$ represents the square of the infinitesimal element of arclength $ds$.
One writes   $ds^2 = \Sigma g_{ij} (q) dq^i dq^j$.  The length
of a path is defined to be $\int ds = \int \sqrt{ \Sigma g_{ij} (q(t) \dot q^i \dot q^j } dt$
and is independent of how the curve $q(t)$ is parameterized.    
The geodesic equations are the EL equations for $L$. 
A judicious application of the  Cauchy Schwartz inequality 
shows that a curve is a minimizer for  the standard problem
of the calculus of variations  for $L$
if and only if it both minimizes the total length among all paths
joining the specified endpoints and is parameterized so as to have constant speed:
$L(q(t), \dot q(t)) = const.$. 
\end{example}

\subsubsection{Heuristic Lagrangian proof of JM formulation of mechanics}  Write $U = -V$ for the negative of the
potential.  Then  Newton's equations
are the Euler-Lagrange equations for the Lagrangian $L(q,v) = K(v) + U (q)$.
At energy $E =0$ we have $U \ge 0$ on the Hill region.
  Using $a^2 + b^2 \ge 2ab$
with equality iff $a =b$ we get that $K +U  \ge 2 \sqrt{K U}$
with equality if and only if $K = U$, which is to say, iff $E = 0$.
Now integrate with respect to $t$ along any path. 
We get $\int L dt \ge \int 2 \sqrt{K U} dt$.  But $2 \sqrt{KU}dt = \sqrt{ 2U} |\dot q| dt = ds_{JM}$ 
is   the zero energy JM arclength for curves  parameterized by Newtonian  time $t$.  So
$$\int_{q[0,T]} L dt \ge \int_{q([0,T])} ds_{JM} $$
  with equality between the integrals if and only if $E(q(t), \dot q(t)) = 0$   along the path
  of integration $q$.  Now   recall  the connection to  the calculus of variations.  The EL equations
are the extremal equations associated
to the action, the integral of $L$ over paths.    In the
formulation above of the standard problem in the calculus of variations
  we fixed both the  endpoints $q_0, q_1$ of paths, {\it and also the time of travel $T$}
between these endpoints.  Relax the time travel condition so as to allow any   time of
flight between the two points.  We have shown that minimizers for this relaxed problem must have zero energy.
They will also  satisfy the EL equations by the standard argument of
the calculus of variation, and, being JM arclength minimizers,  must be JM geodesics.
We have shown that   the free-time action  minimizers must (a) have zero energy and (b) be a zero-energy
JM geodesic.   

``QED''  
 
The astute reader will see various problems implicit in the above ``proof''.  We will not
try to fix them, but leave them as challenges to the reader.    Instead, we will describe
the simple trick which allows us to promote this  heuristic argument to any energy $E$.

If we replace the potential $V$ by $V+ c$,  $c$ a constant, we get the
same Newton's equations.  It follows that if we replace $U$ by $U+E$ for $E$ a constant
in the Lagrangian we still get Newton's equations as the EL equations. 
Now $K = a^2$ but $b = \sqrt{U +E}$ where we need $U + E \ge 0$ so that
we are working on the Hill's region.  We get that the modified action   $\int (K + U + E) dt$
is greater than or equal to the energy $E$ JM arclength $\int 2 \sqrt{K} \sqrt{U + E} dt$,
and that the free-time minimizers of the modified action must be
both solutions to Newton's equation having energy $E$, when parameterized
by time, and a  JM-geodesic with respect to the energy $E$ JM metric.

\subsection{  Hamiltonian formalism.}      In the Lagrangian formulation  we work with positions and velocities. In the Hamiltonian version  we work with positions and momenta.      To explain the difference 
it will help to replace $\R^n$ by an abstract real finite-dimensional vector space  $\V$. 
This $\V$ is meant to encode the position of our ``mechanical system''.   
If $q(t)$ is a curve in
$\V$ its derivative $v(t) = \dot q (t)$ is also a curve in  $\V$.  So velocities are  
vectors $v \in \V$.  Momenta on the other hand  are dual vectors  $p \in \V^*$.         The Lagrangian
is   a function on $L: \V \times \V \to \R$ while the Hamiltonian is a function
 $H: \V \times \V^* \to \R$.   

The Legendre transformation 
 $$FL: \V \times \V \to \V \times \V^*;    FL(q,v) = (q,  \dd{ L}{v} (q,v))$$
mediates between the Lagrangian and Hamiltonian formalisms.    Some words 
are in order regarding   the partial derivative $\dd{ L}{v}$ and why it takes values
in $\V^*$.  
Freeze   $q$ and 
consider the function  $f(v) = L(q, v)$, $f: \V \to \R$.   Then 
$df(v_0) (h) = \frac{d}{d \epsilon} f(v_0 + \epsilon h)$, valid for all $h \in \V$
and so $df(v_0) =  \dd{ L}{v} (q,v_0 ) \in \V^*$
What is important   is that $p =  \dd{ L}{v} (q,v) \in \V^*$.

\begin{exercise} A. Suppose that   $\V = \R$ is one-dimensional. 
and $L(q,v) = \frac{1}{2} m v^2 - V(q)$.  Compute that  $FL (q, v) = (q, m v)$
so that momentum is given by the usual formula $p = m v$.   

B.  Suppose that $\V = \R^n$ and $L(q, v) = \frac{1}{2} \Sigma m_a (v^a)^2 - V(q)$.
Then the  momentum components of $p = \dd{ L}{v} (q,v)$ are given by $p_a = m_a v^a$. Here we have  invoked
the identification $\R^n = (\R^n)^*$ 
induced by the canonical basis, or,
what is the same, by the standard inner product.  In this way   $p = (m_1 v^1, \ldots , m_n v^n) \in \R^n$.

C.  Suppose that $\V$ is endowed with a Euclidean inner product $\langle, \cdot,  \cdot \rangle$
and that $L(q,v) = \frac{1}{2} \langle v, v \rangle - V(q)$.  Show
that   $p = \dd{L}{v} (q,v)$ is given  by metric duality:
$p: \V \to \R$ is the linear functional $p(w) = \langle v, w \rangle$.

D.  Return to the example of Riemannian geometry above.  
Show that $\dd{L}{v}(q,v)_i = \Sigma g_{ij} (q) v^j$.  In other words, if 
I write $\langle v, w \rangle_q = \Sigma g_{ij} (q) v^i w^j$ then
$FL(q,v) = (q, p)$ with $p(w) =  \langle v, w \rangle_q$.

\end{exercise}

\begin{exercise}
Show that the EL equations are the ODEs 
$\frac{d}{dt} FL(q,v) = (v, \dd{L}{q}  (q,v))$
on $\V \times \V^*$.  
\end{exercise}

The differential $dV(q)$ of $V$ at a point $q \in \A$ is an element
of the dual space $\V^*$  defined by
$dV(q)(h) = \frac{d}{d \epsilon} V(q + \epsilon h)$.
The gradient of $V$ is  defined by
$$dV(q) (h) = \langle \nabla V (q), h \rangle.$$

To get the Hamiltonian formalism we use the Legendre transformation to change variables from $(q,v)$
to $(q,p)$.  This is most simply done when the transformation is invertible.  
Suppose this is so.    Define the Hamiltonian
 $H: \V \times \V^* \to \R$
  by $H(q, p) = p (v) - L(q,v)$ where $(q,v) = FL^{-1} (q,p)$.    

Assume now that  the Legendre transformation
is invertible.  Then, as the reader can verify without much difficulty,
 the EL equations can be rewritten in the equivalent form
$$\dot q =  \dd{H}{p}(q,p),  \dot p = - \dd{H}{q} (q,p).$$
known as {\it Hamilton's equations}. 
Note the consistency:  $ \dd{H}{q}  (q,p) \in \V^*$ since the partial derivative (differential)
 is with respect to
  $q \in \V$, while $\dd{H}{p} \in \V^{**} = \V$
since the differential is with respect to $p \in \V^*$.

\begin{exercise}  

Return to example C above where 
the Lagrangian is  $L(q, v) = K(v) - V(q)$
with  $K(v) = \frac{1}{2} m \langle v , v \rangle$,
with $q, v \in \V$   a real vector space   $\V$  endowed with a Euclidean structure
$\langle \cdot,  \cdot \rangle$.   The corresponding 
EL equations are equivalent to  the first order form of  Newton's equations  
$$\ddot q = - \nabla V (q)$$
Note the  Euclidean structure
enters in obtaining the gradient of $V$ via $dV(q) (w) = \langle \nabla V(q), w \rangle$. 
The Hamiltonian is $H(q,p) = K^* (p) + V(q)$
where $K^*$ is the quadratic
form associated to the dual metric induced by the inner product
on $\V$ which led to the kinetic energy $K$.
The reader can verify directly that Hamilton's
equations, Newton's equations and the   EL equations
  are rewrites of   each other, the Legendre transformation $p = \langle v, \cdot \rangle$ being
  the change of variables from velocity to momentum.

Show that if $L(q,v) = \frac{1}{2} \Sigma g_{ij} (q) v^i v^j$
is the Lagrangian for a Riemannian geometry then
$H(q,p) = \frac{1}{2} \Sigma g^{ij} (q) p_i p_j$ is
the corresponding Hamiltonian.   Thus the geodesics are the 
configuration projections  $q(t)$ of solutions $(q(t), p(t))$  to this $H$'s Hamilton's equations.
Show that a geodesic is parameterized by arclength if and only
if $H(q(t), p(t)) = 1/2 $.

Consider the special case where $L(q,v) = \frac{1}{2} \lambda (q) |v|^2$
where $|v|^2 = \Sigma (v^i)^2$.  Show that $H(q,p) = \frac{1}{2} \frac{1}{\lambda (q)} |p|^2$.
 
 \end{exercise}

{\sc  A wee bit of  symplectic geometry}.  
The domain  $\E = \V \times \V^*$ of Hamilton's equations
is known as phase space. 
and write vectors of $\E$ in the form   $(V_q,  V_p)$
so that $V_q \in \V,  V_p \in \V^*$.  
\begin{definition}
The canonical two-form on $\E$
is the bilinear  skew-symmetric form 
$\omega:  \E \times \E  \to \R$
given by  $\omega ( (V_q, V_p), (W_q, W_p)) =  (W_p (V_q)  - V_p (W_q)$. 
$\omega ( (q,p), (q',p')) =  p'(q) - p(q')$. 
\end{definition}

The form $\omega$ is non-degenerate, which means
that the linear  map 
$$J:  \E \to \E^* \text{ defined by }  Z \mapsto \omega (Z, \cdot)$$
is a vector space isomorphism. 

\begin{definition}  A symplectic form 
is a non-degenerate skew-symmetric form on a vector space.
A symplectic vector space is a vector space endowed with a
symplectic form.
\end{definition}

A symplectic vector space must have even dimensions.
Any two symplectic vector spaces (like any two inner product spaces)
of the same dimension are isomorphic as symplectic vector spaces.  

If $(\E, \omega)$ is a symplectic vector space and  function $H : \E \to \R$ is
a smooth function then we can define  its
Hamiltonian vector field, or ``symplectic gradient''
by
$$\omega (sgrad H,  Z) = dH (z) (Z).$$
valid for all vectors $Z$ in $\E$ and all points $z \in \E$.
Here $dH(z) \in \E^*$ is the differential of $H$.  
Equivalently
$$sgrad(H) = J dH$$

\begin{exercise}
Verify that $sgrad H = (\dd{H}{p},  - \dd{H}{q})$

Verify that $dH(sgrad H) = 0$.
\end{exercise}

\begin{definition}  A level set $\{F = c \}$ of a smooth
function on a vector space is said to be ``non-degenerate''
if $dF(p) \ne 0$ whenever $F(p) = c$.
\end{definition}

The implicit function theorem asserts
that a non-degenerate level set  $\Sigma = \{p: F(p) = c \}$ is a smooth
hypersurface and that the tangent spaces
to $\Sigma$ are swept out by the kernels of $dF$.
In other words, $T_p \Sigma = ker(dF(p))$
whenever $p \in \Sigma$.

\begin{cor}  If two different functions on phase space
share the same non-degenerate level set, then their Hamiltonian
vector fields are parallel along this level set.
\end{cor}

Proof.  Write $H$ and $F$ for the two functions and
$\Sigma$ for their common level set.
Thus  $\Sigma = \{H = c_1 \} = \{F = c_2\}$
for constants $c_1, c_2$.  The tangent space
to $\Sigma$ at any point $z$ is a codimension
$1$ hyperplane equal to both the kernel of $dH(z)$
and to the kernel of $dF(z)$.  It follows
that $dF(z) = \lambda (z) dH (z)$ holds along $\Sigma$
for some smooth nowhere zero scalar function $\lambda: \Sigma \to \R$.
Applying $J$ we get that
$sgrad G (z) = \lambda (z) sgrad H(z)$,
valid along $\Sigma$.  QED

\subsection{ Hamiltonian proof of the JM reformulation of Mechanics}

The geodesics for the JM metric at energy $E$ are generated
by the Lagrangian $L_{JM} = \frac{1}{2} (2 (E- V(q))) |v|^2$.
Its corresponding Hamiltonian is $F(q, p) = \frac{1}{2} \frac{1}{2 (E-V(q))}|p|^2$.
To get geodesics parameterized by arclength  we must set $F = \frac{1}{2}$.
Note:  $F = \frac{1}{2} \iff \frac{1}{2 (E - V(q)} |p|^2 = 1 \iff \frac{1}{2} |p|^2 = E - V(q)
\iff \frac{1}{2} |p|^2  + V(q)  = E$.  We have shown that the level
set $F = 1/2$ equals the level set $H = E$.  By the above corollary, this
proves that $sgrad (F)(z)  = \lambda (z)  sgrad (H)(z)$ holds
for $z \in \Sigma$,  assuming that the level set is non-degenerate. 
Indeed, a peek into the proof of the corollary shows that we can relax
the assumption that the entire level set be non-degenerate so that the corollary
holds at all points $z$ for which both $dF(z)$ and $dH(z)$ are nonzero.
Since both functions are quadratic in $p$ this non-degeneracy holds
wherever $p \ne 0$, which is to say, away from the Hill boundary. 
(The level set is non-degenerate if and only there are no critical points
of the potential for which $V(q) =E$.)     

To get the   reparameterization formula relating $s$ and $t$ we find $\lambda (q)$.
And to find $\lambda (q)$ we write out the q-part of the two Hamilton's  equations.
We have, for Newtonian version that  $\dot q = p$ (our masses are all ``$1$'').  
On the other hand, for the JM geodesic equations we have that $\frac{dq}{ds} = \frac{1}{2(E- V(q)} p$.
I claim that it follows that $\lambda (q) = \frac{dt}{ds} =  \frac{1}{2(E- V(q)}$.
Indeed, we can write  Newton's equations as $(\frac{dq}{dt}, \frac{dp}{dt})  = (X_1, X_2)$
where  $sgrad(H) = (X_1, X_2)$  and we can write the   geodesic equations as  $(\frac{dq}{ds}, \frac{dp}{ds})  = (Y_1, Y_2)$
where  $sgrad (F) = (Y_1, Y_2)$.  By the corollary,  along $\Sigma$ we have that $Y_i = \lambda (q) X_i$.
But we have just shown that $X_1 = p$ while $Y_1 =  \frac{1}{2(E- V(q)} p$.
Now use $\frac{d}{ds} = \frac{dt}{ds} \frac{d}{dt}$.

 QED

 \section{Acknowledgements}  I  thank Alain Albouy and Andreas Knauf
 for  conversations and history and Gil Bor with help making the pictures.    I  thank the support of 
 a Simons  Foundation Travel grant, award number  674861,  which helped to make this work possible.

 \end{document}